\documentclass[12pt]{article}
\usepackage{amsmath,amsfonts,amsthm,amscd}
\newtheorem{theorem}{Theorem}

\newtheorem{proposition}{Proposition}
\newtheorem{lemma}{Lemma}
\newtheorem{definition}{Definition}
\newtheorem{corollary}{Corollary}

\newcommand{\sk}{\sigma_k^{1/k}(S^n)}
\newcommand{\p}{\varphi}
\numberwithin{equation}{section}
\numberwithin{theorem}{section}
\numberwithin{proposition}{section}
\numberwithin{lemma}{section}
\numberwithin{claim}{section}
\numberwithin{corollary}{section}
\begin{document}
\title{Volume comparison and the $\sigma_k$-Yamabe problem}
\author{ {Matthew J. Gursky}\thanks{Supported in part by NSF 
Grant DMS-0200646} \and 
{Jeff A. Viaclovsky}\thanks{Supported in
part by NSF Grant DMS-0202477.}}
\date{}
\maketitle
\begin{abstract}
In this paper we study the problem of finding a conformal metric with the property that
the $k$-th elementary 
symmetric polynomial of the eigenvalues of its Weyl-Schouten tensor is 
constant.  A new conformal invariant involving maximal 
volumes is defined, and this invariant is then used in several cases 
to prove existence of a solution, and compactness of the space 
of solutions (provided the conformal class admits 
an {\em{admissible}} metric). In particular, the problem is 
completely solved in dimension four, and in dimension 
three if the manifold is not simply connected. 
\end{abstract} 

\section{Introduction}

Let $(M^n,g)$ be a smooth, closed Riemannian manifold of dimension $n$.  We denote the
Riemannian curvature tensor by $Riem$, the Ricci tensor by $Ric$, and the
scalar curvature by $R$.  In addition, the {\it Weyl-Schouten tensor} is defined
by 
\begin{align}
\label{WStensor}
A = \frac{1}{(n-2)} \Big( Ric - \frac{1}{2(n-1)} R g \Big).
\end{align}
Note that under the action of $O(n)$ the curvature tensor can be decomposed as
\begin{align}
\label{rot}
Riem = W + A \odot g,
\end{align}
where $W$ denotes the Weyl curvature tensor, 
and $\odot$ the Kulkarni-Nomizu product \cite {Besse}.  Since the Weyl tensor is
conformally invariant, an important consequence of the 
decomposition (\ref{rot}) is that the tranformation of the Riemannian curvature tensor
under conformal deformations of metric is completely determined by the transformation
of the symmetric $(0,2)$-tensor $A$.  

In \cite{Jeff1}, the second author
initiated the study of the fully nonlinear equations arising from
the transformation of $A$ under conformal deformations.  More precisely, let $g_{u}
=e^{-2u}g$ denote a conformal metric, and consider the equation 
\begin{align}
\label{sigmak}
\sigma_k^{1/k}(g_u^{-1}A_u) = f(x),
\end{align}
where $\sigma_k : \mathbf{R}^n \to \mathbf{R}$ denotes the elementary symmetric polynomial
of degree $k$, $A_u$ denotes the Weyl-Schouten tensor with respect to the metric 
$g_{u}$, and $\sigma_k^{1/k}(g_u^{-1}A_u)$ means $\sigma_k(\cdot)$ applied to the eigenvalues of 
the $(1,1)$-tensor $g_u^{-1}A_u$ obtained by "raising an index" of $A_u$.

To simplify our formulas we usually interpret $A_u$ as a bilinear form on the 
tangent space with inner product $g$ (instead of $g_u$).
That is, once we fix a background metric $g$, 
$\sigma_k(A_u)$ means $\sigma_k(\cdot)$ applied to the eigenvalues of the $(1,1)$-tensor 
$g^{-1}A_u$.  
To understand the practical effect of this convention, recall that
$A_u$ is related to $A$ by the formula
\begin{align}
A_u = A + \nabla^2 u + du \otimes du
- \frac{1}{2}| \nabla u|^2 g
\label{Achange}
\end{align} 
(see \cite{Jeff1}).  Consequently, (\ref{sigmak}) is equivalent to 
\begin{align}
\label{hessk}
\sigma_k^{1/k}(A + \nabla^2 u + du \otimes du
- \frac{1}{2}| \nabla u|^2 g) = f(x)e^{-2u}.
\end{align}
Note that when $k=1$, then $\sigma_1(g^{-1}A) = trace(A) = \frac{1}{2(n-1)}R$.  Therefore, 
(\ref{hessk}) is the
classical problem of prescribing scalar curvature.   This equation is semilinear
elliptic; however, when $k > 1$ equation (\ref{hessk}) is fully nonlinear but not necessarily
elliptic.

To explain the ellipticity properties of (\ref{hessk}),
following G\aa rding \cite{Garding} and Cafarelli-Nirenberg-Spruck 
\cite{CNSIII} we let 
$\Gamma_k^{+} \subset \mathbf{R}^n$ denote the component of
$\{x \in \mathbf{R}^n | \sigma_k(x) > 0 \}$ containing the positive cone $\{ x \in
\mathbf{R}^n | x_1 > 0,..., x_n > 0    
\}$. 
In terms of the cones $\Gamma_k^{+}$,
ellipticity can be characterized in the following manner 
(see \cite{Jeff1}):  If the eigenvalues of $A = A_g$ are everywhere in 
$\Gamma_k^{+}$, 
and if $u$ is a solution to (\ref{hessk}), then $u$ is an
elliptic solution. 
This fact is a consequence of the convexity of the cones 
$\Gamma_k^{+}$.  
Following the usual practice, we
will say that a metric $g$ is {\it k-admissible}
if the eigenvalues of $A = A_g$ are in $\Gamma_k^{+}$, and
we then write $g \in \Gamma_k^{+}(M^n)$.  

The general problem of solving (\ref{hessk}) with $f(x) = constant$ is referred to as
the {\it $\sigma_k$-Yamabe problem}.  It will be convenient to normalize the value of this
constant, so that the round metric on the sphere is a solution (with no need of 
rescaling): 
\begin{align}
\label{skYamabe}
\sigma_k^{1/k}(A + \nabla^2 u + du \otimes du
- \frac{1}{2}| \nabla u|^2 g) = \sk e^{-2u},
\end{align}
where $\sk = \sigma_k^{1/k}(A_0)$, and $A_0$
is the Weyl-Schouten tensor of the round metric on $S^n$.
We remark that the associated
equation is variational when $k = 1$ or $k = 2$, but in general not when $k > 2$ (see
\cite{Jeff1}).

  The variational approach to the classical Yamabe problem lead to 
the definition of the {\it Yamabe invariant} $Y(M^n,[g])$ of a conformal class of metrics:
\begin{align}
Y(M^n,[g]) \equiv \inf_{\tilde{g}\in[g]}(vol(\tilde{g}))^{-(n-2)/n}
\int R_{\tilde{g}}dvol_{\tilde{g}}.
\end{align}
It is a result of Aubin that $Y(M^n,[g]) \leq Y(S^n,g_0)$, where $g_0$ denotes the 
round metric, and  
when strict inequality holds, existence and compactness of solutions
can be easily established (see \cite{LeeandParker}).
Thus, the resolution of the classical Yamabe problem is equivalent to the result, due in some
cases to Aubin (\cite{Aubin}) and in the remaining cases to Schoen (\cite{Schoen4}), that 
equality holds only when the manifold is conformally equivalent to the sphere.

Our first goal in this paper is to define a new conformal invariant associated to 
equation (\ref{hessk}) when $k \geq  n/2$. 

\begin{definition}
Let $(M^n,g)$ be a compact $n$-dimensional Riemannian manifold. For $n/2 \leq k \leq n$ we define
the {\em $k$-maximal volume} of $[g]$ by   
\begin{align}
\label{Lambdak}
\Lambda_k(M^n,[g]) = \sup \{vol(e^{-2u}g) | e^{-2u}g \in \Gamma_k^{+}(M^n)\mbox{ with }
\sigma_k^{1/k}(g_u^{-1}A_u) \geq \sk \}.
\end{align}
If $[g]$ does not admit a $k$-admissible metric, we set $\Lambda_k(M^n,[g]) = +\infty$.
\end{definition}

By recent work of the second author with P. Guan and G. Wang \cite{GVW}, a $k$-admissible metric $g$ with
$k>n/2$ necessarily has positive Ricci curvature.  In fact, their result is quantative, in the
sense that once we make the normalization $\sigma_k^{1/k}(g^{-1}A_g) \geq \sk$ a (sharp) lower bound
for the Ricci curvature is given (see Section 4).  Using Bishop's inequality, it follows that the invariant
$\Lambda_k$ is non-trivial when $k > n/2$:

\begin{proposition}
\label{Lambdafinite}
If $[g]$ admits a $k$-admissible metric with $k > n/2$, then there is a constant $C = C(n)$ such that
$\Lambda_k(M^n,[g]) < C(n)$.
\end{proposition}

When $k = n/2$ the situation is more complicated.  For example, if $(M^n,g)$ is locally conformally flat ($LCF$)
and $n$ is even, then the integral
\begin{align}
\label{k=n/2}
\int_{M^n} \sigma_{n/2}(g^{-1}A)dvol
\end{align}
is conformally invariant; see \cite{Jeff1}.  Therefore, if $g \in \Gamma_{n/2}^{+}(M^n)$ satisfies 
$\sigma_{n/2}^{2/n}(g^{-1}A) \geq \sigma_{n/2}^{2/n}(S^n)$, then 
\begin{align*}
\int_{M^n} \sigma_{n/2}(g^{-1}A)dvol \geq \sigma_{n/2}(S^n)vol(g).
\end{align*}
Consequently, the maximal volume of $[g]$ is finite.  In fact, we can say more: since 
the assumption of $k-$admissibility with
$k > n/2$ already implies that the Ricci curvature is positive, if
$(M^n,g)$ is $LCF$ then by Kuiper's theorem (\cite{Kuiper}) it must be a 
space form.  Since the dimension is even, by Synge's theorem $(M^n,g)$
is conformally equivalent to $S^n$ or $\mathbf{RP}^n$.  Finally, 
by Proposition 8 in \cite{Jeff1} and the
Chern-Gauss-Bonnet formula it follows
 that $\Lambda_{n/2}(M^n,[g]) = vol(S^n)$ or $\Lambda_{n/2}(M^n,[g]) = 
\frac{1}{2}vol(S^n)$, depending on whether $(M^n,g)$ is conformally equivalent to the
sphere or projective space.

In four dimensions the integral (\ref{k=n/2})
is {\em always} conformally invariant, so the preceding argument can be applied to
show the finiteness of $\Lambda_2(M^4,[g])$ for any conformal four-manifold which admits a $2$-admissible metric
(see Theorem \ref{4dest} below for a sharp version of this result). In 
general, however, it is unclear whether $\Lambda_{n/2}$ is finite. 


  In analogy with the classical Yamabe problem, when our invariant is strictly less than 
the value
obtained by the round metric on the sphere we obtain existence of solutions 
to (\ref{skYamabe}):
\begin{theorem}  
\label{Main} 
Let $(M^n,g)$ be a closed $n$-dimensional Riemannian manifold satisfying 
\begin{align}
\label{subcrit}
\Lambda_k(M^n,[g]) < vol(S^n),
\end{align}
where $vol(S^n)$ denotes the volume of the round sphere.  
Then $[g]$ admits a solution $g_u = e^{-2u}g$ of (\ref{skYamabe}).  
Furthermore, the set of solutions of (\ref{skYamabe}) is compact in the 
$C^{m}$-topology for any $m \geq 0$.
\end{theorem}

Despite the parallels with the Yamabe problem, Theorem \ref{Main} 
can only be considered satisfying if
the condition (\ref{subcrit}) is known to be sharp.  Although we conjecture this to be the case in general,
we can only substantiate it in dimensions three and four.  In each case the techniques for
proving  sharp estimates of $\Lambda_k(M^n,[g])$ are quite different in spirit.  

In three dimensions our estimate follows from the volume comparison theorem of Bray (\cite{Bray}). 
We will give a precise statement of his result later; for now we simply state the
consequence for our invariant.

\begin{theorem}
\label{3dest}
Let $(M^3,g)$ be a closed Riemannian three-manifold, and assume $[g]$ admits a
$k$-admissible metric with $k = 2$ or $3$.  Then
\begin{align}
\label{sharp3}
\Lambda_k(M^3,[g]) \leq vol(S^3).
\end{align}
 
\end{theorem}

The proof of this result allows an important refinement of inequality (\ref{sharp3}).  As a consequence, we
are able to
verify the assumptions of Thorem \ref{Main} whenever $M^3$ is not simply connected:

\begin{theorem}
\label{3dpi1}
Let $(M^3,g)$ be a closed Riemannian three-manifold, and assume $[g]$ admits a
$k$-admissible metric with $k = 2$ or $3$.  Let $\pi_1(M^3)$ denote the fundamental
group of $M^3$.  Then
\begin{align}
\label{sharp3wpi1}
\Lambda_k(M^3,[g]) \leq \frac{vol(S^3)}{\|\pi_1(M^3)\|}.
\end{align}
\end{theorem}

\begin{corollary}
\label{3dxpi1}
Let $(M^3,g)$ be a closed, non-simply connected Riemannian three-manifold.  If $g$ is 
$k$-admissible with $k = 2$ or $3$, then $[g]$ admits a solution 
$g_u = e^{-2u}g$ of (\ref{skYamabe}).  
Furthermore, the set of solutions of (\ref{skYamabe}) 
is compact in the $C^{m}$-topology for any $m \geq 0$.
\end{corollary}

In four dimensions, our estimates of $\Lambda_k$
follow from the sharp integral estimate for $\sigma_2(A)$
due to the first author (\cite{Gursky1}).  

\begin{theorem}
\label{4dest}
Let $(M^4,g)$ be a closed Riemannian four-manifold, and assume $[g]$ admits a
$k$-admissible metric with $2 \leq k \leq 4$.  Then
\begin{align}
\label{sharp4}
\Lambda_k(M^4,[g]) \leq vol(S^4).
\end{align}
Furthermore, equality holds in (\ref{sharp4}) if and only if $(M^4,g)$ is conformally
equivalent to the round sphere.
\end{theorem}

\begin{corollary}
\label{4dcor}
Let $(M^4,g)$ be a closed Riemannian four-manifold, and assume $g$ is a
$k$-admissible metric with $2 \leq k \leq 4$.  Then $[g]$ admits a solution 
$g_u = e^{-2u}g$ of (\ref{skYamabe}). 
Furthermore, if $(M^4,g)$ is not conformally equivalent to the round sphere, then 
the set of solutions of (\ref{skYamabe}) 
is compact in the $C^{m}$-topology for any $m \geq 0$.
\end{corollary}

When $k=2$ the result of Corollary \ref{4dcor} was established in \cite{CGY2}.
Combining Corollary \ref{4dcor} with the four-dimensional
solution of the Yamabe problem \cite{Schoen4}, it follows that the
$\sigma_k$-Yamabe problem is completely solved in four dimensions.
 
Similar to the three-dimensional case,
if we impose certain topological conditions then inequality (\ref{sharp4}) can be
sharpened.  Since the work of Viaclovsky-Guan-Wang cited above shows that a
$k$-admissible metric with $k > n/2$ has positive Ricci curvature, by the classical Bochner
theorem the first de Rham cohomology group $H^1(M^4)=0$.  
On the other hand, if the second de Rham
cohomology group is non-trivial, then the $L^2$-estimates of the Weyl curvature
tensor in \cite{Gursky2} can be used to give sharp estimates of the maximal volume.  To state this
result, let $b^{+}$ ({\it resp.,} $b^{-}$)
 denote the dimension of the largest subspace of
$H^2(M^4)$ on which the intersection form is positive ({\it resp.,} negative) definite.  
Let $\chi(M^4)$ denote the Euler characteristic and $\tau(M^4) = b^{+} - b^{-}$ the
signature of $M^4$.  

\begin{theorem}
\label{bpluscase}
Let $M^4$ be a smooth, compact, orientable four-manifold with $b^{+} > 0$.
If $g \in \Gamma_k^{+}$ with $2 \leq k \leq 4$, then
\begin{align}
\label{sharpbetti}
\Lambda_k(M^4,[g]) \leq \frac{2}{9}\pi^2\left( 2\chi(M^4) + 3\tau(M^4) \right).
\end{align}
In particular, if $2\chi(M^4) + 3\tau(M^4) < 12$, then $\Lambda_k(M^4,[g]) < vol(S^4) = \frac{8}{3}\pi^2$.

Furthermore, equality holds if and only if $[g]$ admits a (positive) K\"ahler-Einstein
metric which attains the maximal volume.  
In this case, $M^4$ is diffeomorphic to either $S^2 \times S^2$, $\mathbf{CP}^2$, or $\mathbf{CP}^2
\# m(-\mathbf{CP}^2)$ with $3\leq m \leq 8$.  
\end{theorem}

In higher dimensions we do not have a sharp estimate of our invariant.  However, the
proof of Theorem \ref{3dpi1} can be adapted to give the following result:

\begin{theorem}
\label{higherd}  
There is a number $N$, depending only on $k$ and $n$, 
with the following property: if $M^n$ is 
a closed $n$-dimensional 
manifold whose fundamental group satisfies $\|\pi_1(M^n)\| > N(k,n)$,
then any k-admissible
metric $g$ satisfies $\Lambda_k(M^n,[g]) < vol(S^n)$.
\end{theorem}

There has been a considerable amount of recent activity devoted to the study of 
(\ref{sigmak}) with $k > 1$ (see \cite{CGY1},\cite{CGY2},\cite{CGY3}, \cite{GuanWang2},
\cite{GuanWang1}, \cite{LiLi1},\cite{LiLi2}, \cite{PM}, \cite{Jeff3},\cite{Jeff2}).
With a few notable exceptions, most of these works consider the case where the background
metric is $k$-admissible.

In \cite{Jeff2}, the second author established global
{\it a priori} $C^1$- and $C^2$-estimates for $k$-admissible solutions 
which depend on $C^0$-estimates.  Since (\ref{hessk})
is a convex function of the eigenvalues of $A_u$, the work of Evans and Krylov (\cite{Evans},
\cite{Krylov}) give $C^{2,\alpha}$ bounds once $C^2$ bounds are known. 
Consequently, one
can derive estimates of all orders from classical elliptic regularity, provided $C^0$-
bounds are known.

Subsequently, Guan and 
Wang (\cite{GuanWang1}) proved local versions of these estimates which only depend on 
a lower bound for solutions.  Their estimates will figure prominently in our analysis.  
Recently, Li and Li (\cite{LiLi1}) proved Harnack
estimates for solutions of (\ref{hessk}), and a classification result 
for entire solutions
on $\mathbf{R}^n$.  Their classification result will also be used in the proof of Theorem \ref{Main}

For global estimates the result of (\cite{Jeff2}) is 
optimal: since (\ref{sigmak}) is invariant under the action
of the conformal group, {\it a priori} $C^0$-bounds may fail for the usual reason (i.e.,
the conformal group of the round sphere).  Some results have managed to distinguish the case
of the sphere, thereby giving bounds when the manifold is not conformally equivalent to $S^n$.
For example, \cite{CGY2} proved the existence of solutions to (\ref{hessk}) 
when $k=2$
and $g$ is $2$-admissible, for any function $f(x)$, provided $(M^4,g)$ is not
conformally equivalent to the sphere.  
In \cite{Jeff2} the second author studied
the case $k=n$, and defined another conformal invariant associated to 
admissible metrics.  When this invariant is below
a certain value, one can establish $C^0$-estimates.  Using this fact he proved the
existence of solutions to (\ref{skYamabe}) on a large class
of conformal manifolds. 

When $(M^n,g)$ is locally conformally flat and $k$-admissible, the article
\cite{LiLi2} gives a compactness result for solutions of (\ref{hessk}) for any $k \geq 1$,
assuming $(M^n,g)$ is not conformally equivalent to the sphere.  
Guan and Wang (\cite{GuanWang2}) used a parabolic version of (\ref{skYamabe}) to
prove global existence (in time) of solutions and convergence to 
a solution of (\ref{skYamabe}).   
As we observed above, the assumption of $LCF$ and
$k-$admissibility with
$k \geq n/2$ implies that  $(M^n,g)$ is conformally equivalent to
a space form.

We conclude the introduction with an outline of the paper.   In Section 2
we lay the groundwork for solving (\ref{hessk}) by introducing a one-parameter family
of auxilary equations.  This requires
us to establish various {\it a priori} estimates, which are contained in
Sections 2 and 3.  These estimates allow us to apply the
degree theory for fully nonlinear equations developed by Li (\cite{Yanyan2}) 
to prove the existence of solutions when $\Lambda_k(M^n,[g]) < vol(S^n)$.  
Finally, in Section 4 we prove
some estimates for the conformal invariant $\Lambda_k(M^n,[g])$.  
\medskip

\subsection{Acknowledgements}

The authors are especially grateful to Hugh Bray for bringing to our attention
the volume comparison theorem in his thesis. 
We also benefitted on several occasions from 
discussions with Pengfei Guan and Yanyan Li.  

\section{The auxilary equation: local estimates}

Let $M^n$ be a closed $n$-dimensional manifold, and suppose $g \in \Gamma_k^{+}(M^n)$.  By rescaling,
we assume that $g$ has unit volume.
Consider the equation
\begin{align}
\label{init}
\sigma_k^{1/k}(\lambda_kg + \nabla^2u + du \otimes du
- \frac{1}{2}| \nabla u|^2 g) = \left(\int e^{-(n+1)u}\right)^{\frac{2}{n+1}},
\end{align}
where $\lambda_k$ is given by
\begin{align}
\label{lambdak}
\lambda_k = {n \choose k}^{-1/k}.
\end{align}
This choice of $\lambda_k$ implies $\sigma_k(\lambda_kg) = 1$. Consequently,
$u \equiv 0$ is a solution of (\ref{init}).

\begin{lemma}
\label{unique} 
$u \equiv 0$ is the unique solution of (\ref{init}).
\end{lemma}
\begin{proof}  This follows from the maximum principle, as explained in Proposition 5 of \cite{Jeff2}. 
Suppose $u$ is a solution of (\ref{init}). 
At a point $x_0$ where $u$ attains its maximum,
$\nabla^2 u(x_0)$ is negative semi-definite and $du(x_0)=0$, so (\ref{init}) implies
\begin{align}
\begin{split}
\left(\int e^{-(n+1)u}\right)^{\frac{2}{n+1}} &= \sigma_k^{1/k}(\lambda_kg + \nabla^2u(x_0)) \\
& \leq \sigma_k^{1/k}(\lambda_kg) \\
& = 1.
\end{split}
\end{align}
Applying a similar argument at the minimum of $u$ we
find 
\begin{align}
\left(\int e^{-(n+1)u}\right)^{\frac{2}{n+1}} \geq 1.
\end{align}
Therefore,
\begin{align}
\label{iss1} 
\left(\int e^{-(n+1)u}\right)^{\frac{2}{n+1}} = 1.
\end{align}
By the Newton-Maclaurin inequality,
\begin{align}
\begin{split}
1 &= \left(\int e^{-(n+1)u}\right)^{\frac{2}{n+1}} \\
&= \sigma_k^{1/k}(\lambda_kg + \nabla^2u + du \otimes du
- \frac{1}{2}| \nabla u|^2 g) \\
& \leq  
\frac{1}{n}{{n}\choose{k}}^{1/k}\sigma_1(\lambda_kg + \nabla^2u + du \otimes du
- \frac{1}{2}| \nabla u|^2 g) \\
&= \frac{1}{n}{{n}\choose{k}}^{1/k}(\lambda_k n + \Delta u - \frac{(n-2)}{2}|\nabla u|^2)\\
&= 1 + \frac{1}{n}{{n}\choose{k}}^{1/k}(\Delta u - \frac{(n-2)}{2}|\nabla u|^2).
\end{split}
\end{align}
Then the maximum principle implies $u$ is a constant, and (\ref{iss1}) forces $u \equiv 0$.  
\end{proof}

For the next Lemma, define the operator
\begin{align}
\label{Psi}
\Psi[u] = \sigma_k^{1/k}(\lambda_kg + \nabla^2u + du \otimes du
- \frac{1}{2}| \nabla u|^2 g) - \left(\int e^{-(n+1)u}\right)^{\frac{2}{n+1}}.
\end{align}
By Lemma \ref{unique}, $u_0 \equiv 0$ is the unique solution of
\begin{align}
\label{psi0=0}
\Psi[u_0] = 0.
\end{align}
Let $\mathcal{L}_{u_0}[h] = \frac{d}{ds}\Psi[u_0 + sh]|_{s=0}$ denote the linearization
of $\Psi[\cdot]$ at the solution $u = u_0$. Then
\begin{align}
\label{linearization}
\mathcal{L}_{u_0}[h] = \gamma_{k,n}\Delta h + 2 \int h,
\end{align}
where $\gamma_{k,n} = (n\lambda_k)^{-1}$.

\begin{lemma}
\label{non-degenerate}
$\mathcal{L}_{u_0}: C^{2,\alpha} \to C^{\alpha}$ is invertible.
\end{lemma}

\begin{proof}
Given $f \in C^{\alpha}$, let $h_1$ be the unique solution of 
\begin{align}
\label{h1}
\gamma_{k,n} \Delta h_1 = f - \bar{f}
\end{align}
satisfying 
\begin{align}
\label{inth1}
\bar{h_1} = 0, 
\end{align}
where bars denote the mean value (recall the background metric has unit
volume).  
If we take $h = h_1 + \frac{1}{2} \bar{f}$, then by (\ref{h1}) and (\ref{inth1})
\begin{align*}
\begin{split}
\mathcal{L}_{u_0}[h] &= \gamma_{k,n}\Delta h + 2\int h \\
&= \gamma_{k,n}\Delta h_1 + 2\int \left( h_1 + \frac{1}{2} \bar{f} \right) \\
&= f - \bar{f} + \bar{f} \\
&= f.
\end{split}
\end{align*}
Using the maximum priciple, it is easy to see that $h$ is in fact the unique solution of 
$\mathcal{L}_{u_0}[h] = f$.  
\end{proof}

We now introduce a one-parameter family of equations connecting equation (\ref{skYamabe})
with equation (\ref{init}).
For $t \in [0,1]$, consider
\begin{align}
\label{eqnt}
\begin{split}
\sigma_k^{1/k}&\left(\lambda_k(1-\psi(t))g + \psi(t)A +  \nabla^2 u + du \otimes du
- \frac{1}{2}| \nabla u|^2 g\right) \\
&= (1-t)\left(\int e^{-(n+1)u}\right)^{\frac{2}{n+1}} + \psi(t) \sk e^{-2u},
\end{split} 
\end{align} 
where $\psi(t) \in C^1[0,1]$ satisfies $0 \leq \psi(t) \leq 1, \psi(0) = 0$,
and $\psi(t) \equiv 1$ for $t \geq \frac{1}{2}$.
From the properties of $\psi(t)$ we see that if $u$ is a solution of (\ref{eqnt}) with $t \geq 
\frac{1}{2}$, then $\sigma_k^{1/k}(A_u) \geq \sk e^{-2u}$.  Therefore,
\begin{align}
\label{bound}
\Lambda_k(M^n,[g]) \geq \sup\{vol(g_u) | u \mbox{ satisfies (\ref{eqnt}) with } t \geq \frac{1}{2}\}.
\end{align}

Since (\ref{eqnt}) admits a unique solution when $t = 0$, we would like to use a degree theoretic
argument to show that it also admits a solution when $t = 1$.  The degree theory developed by Li
(\cite{Yanyan2}) for second order fully nonlinear equations provides a framework for this approach.
We will explain the details in Section 3, but it may help the reader to appreciate
the estimates of this section if we first provide an overview of our plan.

The first step is to compute the Leray-Schauder degree of the solution $u \equiv 0$ of (\ref{init}). 
By Lemmas \ref{unique} and \ref{non-degenerate} this degree is non-zero.
The next step is to appeal to the homotopy invariance of the degree to conclude that (\ref{eqnt})
has a solution when $t=1$.  To justify this, however, we need to establish {\it a priori}
bounds for solutions of (\ref{eqnt}).  
As we shall see, when $t < 1$ the integral term in (\ref{eqnt}) imposes $L^{\infty}$-bounds on 
solutions.  By the {\it a priori} $C^1$- and $C^2$-estimates
of \cite{Jeff2}, along with the aforementioned results of Krylov \cite{Krylov} 
and Evans \cite{Evans}, 
such $L^{\infty}$-bounds will imply bounds on derivatives of all orders.

The conformal invariance of equation (\ref{eqnt}) when $t=1$ leads to predictable
difficulties when deriving estimates with $t$ close to $1$.  As $t \to 1$, we
need to use a standard blow-up procedure in order to show that the 
assumption $\Lambda_k(M^n,[g]) < vol(S^n)$ imposes $L^{\infty}$-bounds on solutions.
The classification of solutions of (\ref{skYamabe}) on Euclidean space 
Li and Li (\cite{LiLi2}) will be important in this respect.  

With this overview in mind, we begin with a basic estimate for solutions of (\ref{eqnt}) with
$t < 1$.

\begin{theorem}
\label{mainest}
For any fixed $0 < \delta < 1$, there is a constant $C = C(\delta,g)$ such that any solution of
(\ref{eqnt}) with $t \in [0, 1-\delta]$ satisfies 
\begin{align}
\label{1stest}
\|u\|_{C^{4,\alpha}} \leq C.
\end{align}
\end{theorem}
\begin{proof}  The proof of this estimate is divided into a few intermediate steps, starting
with an estimate on the minimum of solutions.

\begin{proposition}
\label{prop3.1}
If $u$ is a solution of (\ref{eqnt}) with $t \in [0,1-\delta]$, then there is a
constant $C= C(\delta,g)$ such that
\begin{align}
\label{uinf}
\min_{M^n} u \geq C.
\end{align}
\end{proposition}
\begin{proof}  This Proposition is essentially a corollary of the $\epsilon$-regularity result 
for solutions of (\ref{sigmak}) due to 
Guan and Wang (\cite{GuanWang1}).  However, (\ref{sigmak}) and (\ref{eqnt}) differ by
a constant term; thus we need to clarify some estimates to show that their
argument still works.

We begin by noting that the integral in (\ref{eqnt})
is uniformly bounded for $t \leq 1-\delta$.

\begin{lemma}
\label{intOK}
Let $u$ be a solution of (\ref{eqnt}) with $t \in [0,1)$.  Then there is a
constant $C = C(g)$ such that 
\begin{align}
\label{intterm}
(1-t) \left(\int e^{-(n+1)u}\right)^{\frac{2}{n+1}} \leq C.
\end{align}
\end{lemma}
\begin{proof}   
To see this we apply the maximum principle once again:  At a point 
$x_0$ where $u$ attains its maximum, $\nabla^2 u(x_0)$ is negative semi-definite
and $du(x_0) = 0$, so
(\ref{eqnt}) implies
\begin{align*}
\begin{split}
(1-t)\left(\int e^{-(n+1)u}\right)^{\frac{2}{n+1}} & \leq
(1-t)\left(\int e^{-(n+1)u}\right)^{\frac{2}{n+1}} + \psi(t) \sk e^{-2u(x_0)}\\
&= \sigma_k^{1/k}\left(\lambda_k(1-\psi(t))g + \psi(t)A(x_0) + \nabla^2 u(x_0)\right) \\
&\leq  \sigma_k^{1/k}\left(\lambda_k(1-\psi(t))g + \psi(t)A(x_0)\right) \\
&\leq C.
\end{split}
\end{align*}
This proves the Lemma.
\end{proof} 

\begin{corollary}
Let $u$ be a solution of (\ref{eqnt}) with $t \leq 1 - \delta$.  Then there is a
constant $C = C(\delta,g)$ such that 
\begin{align}
\label{intbound}
\int e^{-(n+1)u} \leq C.
\end{align}
\end{corollary}

We now turn to the proof of (\ref{uinf}), arguing by contradiction.
Suppose to the contrary we have a sequence $\{u_j\}$
of solutions of (\ref{eqnt}) with $t = t_j \leq 1 - \delta$, and that $\min u_j \to 
-\infty$.  At a point $z_j$ where $u_j$ attains its minimum let $exp_{z_j} : B(0,\iota_0/2)
\subset T_{z_j}M^n \approx \mathbf{R}^n \to M^n$ denote the exponential map, where $\iota_0$
is the injectivity radius of $(M^n,g)$.  Let $\epsilon_j$ satisfy $\log \epsilon_j = \min u_j
= u_j(z_j)$, and define 
\begin{align}
\begin{split}
T_j(x) &= exp_{z_j}(\epsilon_jx), \\
g_j &= \epsilon_j^{-2}T_j^{*}g, \\
\tilde{u}_j(x) &= (T_j^{*}u_j)(x) - \log\epsilon_j \\
&= u_j(exp_{z_j}(\epsilon_jx)) - \log\epsilon_j.
\end{split}
\end{align}
Then each $\tilde{u}_j$ is defined on $B(0,\epsilon_j^{-1}\iota_0/2) \subset \mathbf{R}^n$ 
and satisfies $\tilde{u}_j(x) \geq 0, \tilde{u}_j(0) = 0$.  In addition, since $u_j$ satisfies (\ref{eqnt}),
$\tilde{u}_j$ satisfies
\begin{align}
\label{euceqn}
\begin{split}
\sigma_k^{1/k}&\left((\lambda_k(1-\psi(t_j))g_j + \psi(t_j)A_j) 
 + \nabla^2\tilde{u}_j + d\tilde{u}_j \otimes d\tilde{u}_j
 - \frac{1}{2}|\nabla \tilde{u}_j|^2g_j\right)\\
&= \epsilon_j^2 (1-t_j)\left(\int e^{-(n+1)u_j}\right)^{\frac{2}{n+1}} + \psi(t_j) \sk e^{-2\tilde{u}_j},
\end{split}
\end{align}
where $A_j = A_{g_j}$, 
and the covariant derivatives in (\ref{euceqn}) are with respect to
$g_j$.  Note that $g_j$ converges to the Euclidean metric $ds^2$ on compact sets in the $C^m$-
topology, for any $m\geq1$.    

Next we claim that for any $\rho > 1$,
there is a constant $C = C(\rho,g)$ such that 
\begin{align}
\label{tildeuc1}
\max_{B(0,\rho)}|\nabla \tilde{u}_j|^2 \leq C.
\end{align}
This estimate is a consequence of the local $C^1$-estimate of Guan and Wang:

\begin{lemma}
(See \cite{GuanWang1}, Proposition 2)  Let $u \in C^3$ be an admissible solution of 
\begin{align}
\label{gweqn}
F(u) = \sigma_k^{1/k}(A + \nabla^2 u + du \otimes du
- \frac{1}{2}| \nabla u|^2 g) = f(x)e^{-2u}
\end{align}
on $B(0,2\rho)$, where $\rho > 0$.  Then there is a constant
$C(k,n,\rho,\|g\|_{C^3(B(0,\rho))},\|f\|_{C^3(B(0,\rho))})$
such that 
\begin{align}
\label{gwc1}
|\nabla u|^2(x) \leq C(1+e^{-2\inf_{B(0,\rho)}u})
\end{align}
for all $x \in B(0,\rho/2)$.  
\end{lemma}

In our case, $\tilde{u}_j$ satisfies
\begin{align}
F(\tilde{u}_j) = \epsilon_j^2 (1-t_j)\left(\int e^{-(n+1)u_j}\right)^{\frac{2}{n+1}} + \psi(t_j) \sk e^{-2\tilde{u}_j}.
\end{align}
If we imitate the proof of \cite{GuanWang1}, the only necessary changes
appear in the estimates of inequality (13) of Proposition 2 in 
\cite{GuanWang1}.  More specifically,
Guan and Wang estimate the term
\begin{align}
\begin{split}
\label{gwgw}
\sum_{l} F_l u_l &= \sum_{l} e^{-2u}(f_lu_l - 2fu_l^2)\\
&\geq -C(1+e^{-2u})|\nabla u|^2,
\end{split}
\end{align}
where the subscript $l$ denotes $\frac{\partial}{\partial x_l}$.   Since our definition of $F$
differs only by a constant term, we can literally copy their argument to obtain the same estimate
for $\tilde{u}_j$:
\begin{align}
\max_{B(0,\rho)}|\nabla \tilde{u}_j|^2 \leq C(\rho,g,\min \tilde{u}_j).
\end{align}
Of course, in our case $\tilde{u}_j \geq 0$, and so (\ref{tildeuc1}) follows.

Combining the gradient bound (\ref{tildeuc1}) with the condition $\tilde{u}_j(0) = 0$ we
see that 
\begin{align}
\label{lower}
\min_{B(0,1)} e^{\tilde{u}_j}  \geq C(g) > 0.
\end{align}
On the other hand, pulling back to $M^n$ by $T_j^{-1}$ and using the
integral bound (\ref{intbound}) we have
\begin{align*}
\begin{split}
\int_{B(0,1)} e^{-(n+1)\tilde{u}_j}dvol_{g_j}  &= \epsilon_j
\int_{B(z_j,\epsilon_j)} e^{-(n+1)u_j} dvol_g \\ 
& \to 0
\end{split}
\end{align*}
as $j \to \infty$.
Since this contradicts (\ref{lower}), we see that the sequence $\{u_j\}$ must be 
bounded from below.
\end{proof}

\begin{proposition}
If $u$ is a solution of (\ref{eqnt}) with $t \leq 1-\delta$, then
there is a constant $C=C(\delta,g)$ such that
\begin{align}
\label{supu}
\max_{M^n} u \leq C.
\end{align}
\end{proposition}
\begin{proof}
As we explained in the proof of Theorem \ref{mainest}, the localized gradient 
estimate of Guan and
Wang can be adapted to equation (\ref{eqnt}), giving the bound
\begin{align}
\label{gradu}
\max_{M^n} |\nabla u| \leq C(1 + e^{-2\min u}) \leq C(\delta,g).
\end{align}
This immediately implies the Harnack inequality
\begin{align}
\label{harnack}
\max_{M^n} u \leq \min_{M^n} u + C.
\end{align}
The upper bound (\ref{supu}) will be a consequence of the following Lemma:

\begin{lemma}
\label{minab}
If $u$ is a solution of (\ref{eqnt}) with $t \in [0,1]$, then there is
a constant $C = C(g)$ such that 
\begin{align}
\label{minabove}
\min_{M^n} u \leq C.
\end{align}
\end{lemma}
\begin{proof}  
Let $x_0$ be a point at which $u$ atttains its minimum.  Then
\begin{align}
\label{1stp}
(1-t)\left(\int e^{-(n+1)u}\right)^{\frac{2}{n+1}} + \psi(t) \sk e^{-2u(x_0)}
\leq \left( 1 + \sk \right)e^{-2\min u}.
\end{align}
At $x_0$, $\nabla^2 u(x_0)$ is positive semi-definite
and $du(x_0) = 0$.  Therefore, 
\begin{align}
\label{2ndstp}
\begin{split}
(1-t)\left(\int e^{-(n+1)u}\right)^{\frac{2}{n+1}} + &\psi(t) \sk e^{-2u(x_0)} \\
&= \sigma_k^{1/k}\left(\lambda_k(1-\psi(t))g + \psi(t)A(x_0) + \nabla^2 u(x_0)\right)\\
&\geq \sigma_k^{1/k}\left(\lambda_k(1-\psi(t))g + \psi(t)A(x_0)\right).
\end{split}
\end{align}
Since $\sigma_k : \Gamma_k^{+} \to \mathbf{R}$ is a concave function (see \cite{Jeff2}, Proposition 1),
\begin{align}
\label{beloww}
\begin{split}
\sigma_k^{1/k}\left(\lambda_k(1-\psi(t))g + \psi(t)A(x_0)\right) &\geq \sigma_k^{1/k}\left(\lambda_k(1-\psi(t))g\right) 
+ \sigma_k^{1/k}\left(\psi(t)A(x_0)\right) \\
&= (1-\psi(t)) + \psi(t)\sigma_k^{1/k}(A(x_0)) \\
&\geq C(g) > 0.
\end{split}
\end{align}
Combining (\ref{1stp}), (\ref{2ndstp}) and (\ref{beloww}) we find
\begin{align*} 
e^{-2\min u} \geq C(g) > 0,
\end{align*}
which implies (\ref{minabove}).
\end{proof}

To complete the proof of Theorem \ref{mainest} we appeal to the 
global {\it a priori} estimates of \cite{Jeff2} (see Propositions
6 and 8):  If $u$ is a solution of (\ref{eqnt}) with $0
\leq t \leq 1 - \delta$, then
\begin{align}
\label{c2t}
\begin{split}
\|\nabla u\|_{\infty} + \|\nabla^2 u\|_{\infty} & \leq C(\|u\|_{\infty}) \\
& \leq C(\delta,g).
\end{split}
\end{align}
As explained in the introduction,  
the work of Evans (\cite{Evans}) and Krylov (\cite{Krylov}) now
give bounds on the Holder norms of the second derivatives of $u$.  Hence, the estimate 
(\ref{1stest}) follows from classical elliptic regularity.
\end{proof}
\end{proof}

\section{Global estimates and existence}

Having established estimates for solutions of (\ref{eqnt}) when $t$ is bounded away from $1$,
we now study what happens as $t \to 1$.  As the title of this section indicates, the
analysis of this case depends on {\it global} invariants of the manifold--namely, $\Lambda_k$--rather
than {\it local} properties of the equation (\ref{eqnt}).

\begin{theorem}
\label{globalest}
Suppose $\Lambda_k(M^n,[g]) < vol(S^n)$.
If $u$ is a solution of (\ref{eqnt}) with $t \in [0,1]$, then there
is a constant $C = C(g)$ such that 
\begin{align}
\label{C4est}
\| u \|_{C^{4,\alpha}} \leq C.
\end{align}
\end{theorem}
\begin{proof}  Like the proof of Theorem \ref{mainest}, we use a blow-up argument.  However,
since the integral bound (\ref{intbound}) degenerates as $t \to 1$ we can no longer
rely on an $\epsilon$-regularity result.  This is to be  expected, given the phenomenon of
bubbling.  In any case, we still begin with an estimate of the lower bound of $u$.

\begin{proposition}
There is constant $C = C(g)$ such that
\begin{align}
\label{globalblower}
\min u \geq -C.
\end{align}
\end{proposition}
\begin{proof}
Once again, we argue by contradiction: 
Suppose to the contrary we have a sequence $\{u_j\}$
of solutions of (\ref{eqnt}) with $t = t_j \to 1$, and that $\min u_j \to 
-\infty$.  At a point $z_j$ where $u_j$ attains its minimum let $exp_{z_j} : B(0,\iota_0/2)
\subset T_{z_j}M^n \approx \mathbf{R}^n \to M^n$ denote the exponential map, where $\iota_0$
is the injectivity radius of $(M^n,g)$.  As before, let $\epsilon_j$ satisfy $\log \epsilon_j = \min u_j
= u_j(z_j)$, and define 
\begin{align}
\begin{split}
T_j(x) &= exp_{z_j}(\epsilon_jx), \\
g_j &= \epsilon_j^{-2}T_j^{*}g, \\
\tilde{u}_j(x) &= (T_j^{*}u_j)(x) - \log\epsilon_j \\
&= u_j(exp_{z_j}(\epsilon_jx)) - \log\epsilon_j.
\end{split}
\end{align}
Then each $\tilde{u}_j$ is defined on $B(0,\epsilon_j^{-1}\iota_0/2) \subset \mathbf{R}^n$ 
and satisfies $\tilde{u}_j(x) \geq 0, \tilde{u}_j(0) = 0$.  In addition, by (\ref{eqnt})
$\tilde{u}_j$ satisfies (\ref{euceqn}):
\begin{align}
\begin{split}
\sigma_k^{1/k}&\left((\lambda_k(1-\psi(t_j))g_j + \psi(t_j)A_j) 
 + \nabla^2\tilde{u}_j + d\tilde{u}_j \otimes d\tilde{u}_j
 - \frac{1}{2}|\nabla \tilde{u}_j|^2g_j\right)\\
&= \epsilon_j^2 (1-t_j)\left(\int e^{-(n+1)u_j}\right)^{\frac{2}{n+1}} + \psi(t_j) \sk e^{-2\tilde{u}_j}.
\end{split}
\end{align}
Note that by Lemma \ref{intOK}, as $j \to \infty$ the integral term above goes to zero:
\begin{align*}
\begin{split}
\epsilon_j^2 (1-t_j)\left(\int e^{-(n+1)u_j}\right)^{\frac{2}{n+1}} &\leq C(g)\epsilon_j^2 \\
&\to 0.
\end{split}
\end{align*}

The localized estimate of Guan and Wang (\cite{GuanWang1}) implies that 
for any $\rho > 1$,
there is a constant $C = C(\rho,g)$ such that 
\begin{align}
\max_{B(0,\rho)}|\nabla \tilde{u}_j|^2 \leq C.
\end{align}
Combining this gradient bound with the condition $\tilde{u}_j(0) = 0$ we
see that 
\begin{align}
\max_{B(0,\rho)} (|\tilde{u}_j| +|\nabla \tilde{u}_j| ) \leq C(\rho),
\end{align}
for any $\rho > 1$.  With this estimate we can appeal to the local $C^2$-estimates of Guan and Wang
(\cite{GuanWang1}, Proposition 3).  Once again, our equation is slightly different, but this time (in light
of the $C^1$-estimates for $\tilde{u}_j$) the required modifications are minor.  We will omit the details.
As a consequence, on any ball $B(0,\rho)$,  $\tilde{u}_j$ satisfies 
\begin{align}
\max_{B(0,\rho)} (|\tilde{u}_j| +|\nabla \tilde{u}_j| + |\nabla^2 \tilde{u}_j|) \leq C(\rho).
\end{align}
It follows from the work of Evans and Krylov that one obtains $C^{2,\alpha}$-estimates for
$\tilde{u}_j$ on any fixed ball, and consequently 
$\{\tilde{u}_j\}$ converges uniformly in the $C^{2,\alpha}$-topology on compact sets to a solution $u$ of
\begin{align}
\label{standard}
\sigma_k^{1/k}(\nabla^2 u + du \otimes du
- \frac{1}{2}| \nabla u|^2 g) = \sk e^{-2u}.
\end{align}
The aforementioned regularity results imply that $u \in C^{\infty}$.

By the classification result of Li and Li \cite{LiLi2}, all solutions of (\ref{standard}) are obtained
by pulling back the round metric on the sphere (and its images under conformal diffeomorphisms)
via stereographic projection.  In particular, 
\begin{align}
vol(e^{-2u}ds^2) = vol(S^n).
\end{align}

\begin{lemma}
\label{volumeisbigg}
\begin{align}
\label{crit}
\liminf_{j}vol(e^{-2u_j}g) \geq vol(S^n).
\end{align}
\end{lemma}
\begin{proof}
Given $\eta > 0$, fix a large ball $B = B(0,\rho_0) \subset \mathbf{R}^n$ such that 
\begin{align}
\int_{B} e^{-n\tilde{u}_j} dvol_{g_j} > vol(S^n) - \eta
\end{align}
for all $j \geq J$.  Pulling back to $M^n$ by $T_j^{-1}$, we have
\begin{align}
\begin{split}
vol(e^{-2u_j}g) = 
\int e^{-nu_j} dvol_g &\geq \int_{B(z_j,\epsilon_j \rho_0)} e^{-nu_j} dvol_g \\
& = \int_{B} e^{-n\tilde{u}_j} dvol_{g_j} \\
& > vol(S^n) - \eta.
\end{split}
\end{align}
This proves the Lemma.
\end{proof}

By equation (\ref{eqnt}), for $t \geq \frac{1}{2}, u_j$ satisfies 
\begin{align}
\sigma_k^{1/k}(A + \nabla^2u_j + du_j \otimes du_j - \frac{1}{2}|\nabla u_j|^2g)
\geq \sk e^{-2u_j}.
\end{align}
Therefore, $g_j = e^{-2u_j}g$ satisfies
\begin{align}
\sigma_k^{1/k}(g_j^{-1} A_{u_j}) \geq \sk.
\end{align}
From Lemma \ref{volumeisbigg} we conclude $\Lambda_k(M^n,[g]) = vol(S^n)$, which is a contradiction.  Therefore,  
the sequence $\{u_j\}$ is bounded below, as claimed. 
\end{proof}
To complete the proof of the theorem, we may argue exactly as in the proof of Theorem \ref{mainest}.  Namely,
the localized $C^1$-estimate of Guan and Wang together with the lower bound (\ref{globalblower}) implies a gradient
bound for $u$, and consequently the Harnack inequality (\ref{harnack}).  We may then apply Lemma \ref{minab}
to conclude that $u$ has an {\it a priori} upper bound.  Higher order estimates follow, just as we described
at the end of the proof of Theorem \ref{mainest}.
\end{proof}

The preceding blow-up argument can be applied, with only minor modifications, 
to prove the compactness of solutions of (\ref{skYamabe}).  The 
details will be omitted.

To establish existence, we apply the degree theory for fully nonlinear 
equations as developed in \cite{Yanyan2}. 
In Section 2 we showed that the Leray-Schauder degree of a solution of (\ref{eqnt}) at $t=0$ is nonzero. 
We remark that equation (\ref{eqnt}) differs from that considered 
in \cite{Yanyan2} only by the presence of the integral 
term. From the compactness established in Theorem \ref{mainest}, this integral 
term is bounded. Furthermore, the proof in \cite{Yanyan2} 
relies on differentiating the equation. Since the integral 
term is a constant, the definition of degree and proof of 
invariance of degree under homotopy are valid for equation (\ref{eqnt}). 
We conclude that the Leray-Schauder degree at $t=1$ (with respect 
to a sufficiently large ball in $C^{4,\alpha}$) is nonzero, and
consequently there exists a solution at $t=1$.

\section{Sharp estimates for $\Lambda_k$}

In this section we prove various estimates for the conformal invariant $\Lambda_k$. 
We begin by describing some general properties which are independent of the 
dimension, then consider the cases $n=3$ and $n=4$ separately.

A basic tool in many of our results is the {\it Newton-Maclaurin inequality} (see \cite{Hardy}): 
if $(\lambda_1,...,\lambda_n) \in \Gamma_k^{+}$ and $k\geq j$, then
\begin{align*}
{n \choose k}^{-1/k}\sigma_k^{1/k}(\lambda_1,...,\lambda_n) \leq {n \choose j}^{-1/j}\sigma_j^{1/j}
(\lambda_1,...,\lambda_n).
\end{align*}
This implies
\begin{lemma}
If $g \in \Gamma_k^{+}(M^n)$ and $j \leq k$, then 
\begin{align}
\label{kandj}
\Lambda_j(M^n,[g]) \geq \Lambda_k(M^n,[g]).
\end{align}
\end{lemma}

As a conseqeunce of this Lemma, in order to estimate $\Lambda_k$ with $k > n/2$ it typically suffices
to estimate $\Lambda_j$, where $j = [\frac{n}{2}] + 1$.

\medskip

\noindent{\bf The proof of Proposition \ref{Lambdafinite}.}
The proof is based on the
sharp inequality of Guan, Viaclovsky, and Wang (\cite{GVW}):  If 
$g \in \Gamma_k^{+}(M^n)$ with $k > n/2$, 
then the Ricci tensor satisfies 
\begin{align}
\label{Riccilower} 
Ric \geq \frac{(2k-n)}{2n(k-1)}Rg,
\end{align}
where $R$ is the scalar curvature of $g$.  
The finiteness of $\Lambda_k(M^n,[g])$ will follow from a lower bound 
for the scalar curvature and the Bishop Comparison Theorem, as 
we now explain.  
 
First, by the Newton-MacLaurin inequality 
\begin{align}
\label{Rbelow}
\sigma_k^{1/k}(g^{-1}A) \leq  \frac{1}{n}{{n}\choose{k}}^{1/k}\sigma_1(g^{-1}A) = c(k,n) R.
\end{align}
If $g$ satisfies 
\begin{align}
\label{slower}
\sigma_k^{1/k}(g^{-1}A) \geq \sk,
\end{align}
then combining (\ref{Rbelow}) and (\ref{slower}) we have
\begin{align*}
R \geq c(k,n)^{-1}\sk > 0.
\end{align*}
Substituting this into (\ref{Riccilower}) gives
\begin{align*}
Ric \geq \frac{(2k-n)}{2n(k-1)}c(k,n)^{-1}\sk g.
\end{align*}
Since $n/2 < k \leq n$, we obtain a lower bound for $Ric$ which only depends on
the dimension:
\begin{align*}
Ric \geq c(n)g.
\end{align*}
By Myer's theorem,
the diameter of $g$ is bounded by a constant $C=C(n)$:   
\begin{align*}
diam(M^n,g) \leq C.
\end{align*}
In addition, by the Bishop comparison theorem the positivity of the Ricci curvature implies the volume of a geodesic ball
of radius $\rho$ in $g$ is bounded by $vol(B(\rho)) \leq C_n \rho^n$.  This fact,
combined with the diameter estimate above, implies that $vol(M^n,g) \leq C(n)$.  Thus,
\begin{align*}
\Lambda_k(M^n,[g]) \leq C(n).
\end{align*}
This completes the proof.

\subsection{n = 3}

We now turn to three dimensions, where the sharp estimates of $\Lambda_k$ are based on the 
following result of H. Bray (\cite{Bray}):

\begin{theorem}[Bray's Football Theorem]
\label{braythesis}
Let $(S^3,g_0)$ be the constant curvature metric on $S^3$ with scalar curvature $R_0$, 
Ricci tensor $Ric_0g_0$, and volume $V_0$.  If $\epsilon \in (0,1]$
and $(M^3,g)$ is any complete smooth Riemannian manifold of volume $V$ satisfying
\begin{align}
\label{Braynormal}
R(g) \geq R_0,
\end{align}
\begin{align}
\label{brayRc}
Ric(g) &\geq \epsilon Ric_0g,
\end{align}
then
\begin{align}
\label{volest}
V \leq \alpha(\epsilon)V_0, 
\end{align}
where
$$
\alpha(\epsilon) = \sup_{\frac{4\pi}{3-2\epsilon}\leq z \leq 4\pi} \frac{1}{\pi^2}\left(
\begin{array}{rr} 
\int_0^{y(z)}(36\pi-27(1-\epsilon)y(z)^{\frac{2}{3}}-9\epsilon x^{\frac{2}{3}})^{-\frac{1}{2}}dx \\
+ \int_{y(z)}^{z^{\frac{3}{2}}}(36\pi-18(1-\epsilon)y(z)x^{-\frac{1}{3}}-9x^{\frac{2}{3}})^{-\frac{1}{2}}dx\end{array}\right),
$$
where 
\begin{align*}
y(z) = \frac{z^{\frac{1}{2}}(4\pi-z)}{2(1-\epsilon)}.
\end{align*}
Furthermore, this expression for $\alpha(\epsilon)$ is sharp.
\end{theorem}

When $\epsilon = 1$, the lower bound on the
scalar curvature (\ref{Braynormal}) follows from the lower bound on the Ricci curvature (\ref{brayRc}), and the result
is equivalent to Bishop's inequality.  Now define
\begin{align*}
\epsilon_0 = \inf \{ \epsilon \in (0,1] | \alpha(\epsilon)=1 \}.
\end{align*}
Bray's theorem is remarkable precisely because $\epsilon_0 < 1$.  Although Bray claimed this fact in his thesis, 
he did not include the proof.  However, he did 
provide compelling numerical evidence suggesting $\epsilon_0 = 0.134...$  This value of 
$\epsilon_0$ corresponds to a rotationally symmetric manifold resembling a football; thus the name.  In any case,
there are currently no rigorous estimates of $\epsilon_0$ from above.

For our purposes we need to know 
that $\epsilon_0 \leq 0.5$.   
To see why, suppose $g \in \Gamma_2^{+}(M^3)$ satisfies 
\begin{align}
\sigma_2^{1/2}(g^{-1}A) \geq \sigma_2^{1/2}(g_0^{-1}A_0).
\end{align}
Then the Newton-Maclaurin inequality implies 
\begin{align}
\label{Rlow0}
R \geq R_0.
\end{align}
In addition, by inequality (\ref{Riccilower}),
\begin{align}
\label{Rclow0}
\begin{split}
Ric(g) &\geq  \frac{1}{6}Rg \\
&\geq \frac{1}{6}R_0g \\
&= \frac{1}{2}Ric_0g.
\end{split}
\end{align}
Therefore, if $\epsilon_0 \leq \frac{1}{2}$, from Bray's theorem we would conclude 
\begin{align}
\label{braysays}
\sigma_2^{1/2}(g^{-1}A) \geq \sigma_2^{1/2}(g_0^{-1}A_0)
 \Rightarrow   vol(M^3,g) \leq vol(S^3).
\end{align}
Consequently,
\begin{align*}
\Lambda_2(M^3,[g]) \leq vol(S^3),
\end{align*}
and Theorem \ref{3dest} would follow.

A similar argument can be used to prove inequality (\ref{sharp3wpi1}), again provided 
$\epsilon_0 \leq \frac{1}{2}$.  Under the assumptions of 
Theorem \ref{3dpi1}, we know from inequality (\ref{Riccilower}) that $g$ has positive
Ricci curvature.  By Meyer's theorem the fundamental group of $M^3$ is finite.  If we 
let $\tilde{M}^3$ denote the universal cover of $M^3$, then $\tilde{M}^3$ is compact
and the volume of $M^3$ and $\tilde{M}^3$ are related by 
\begin{align}
\label{covervol}
vol(\tilde{M}^3,\tilde{g}) = \|\pi_1(M^3)\|vol(M^3,g),
\end{align}
where $\tilde{g}$ denotes the lift of $g$ to $\tilde{M}^3$.  Applying the volume estimate 
(\ref{braysays}) to the cover $(\tilde{M}^3,\tilde{g})$ and using  (\ref{covervol}), 
we arrive at (\ref{sharp3wpi1}).  A similar argument can be used to prove
Theorem \ref{higherd}.

The main result of this subsection is a rigorous proof of the inequality $\epsilon_0 \leq \frac{1}{2}$.
Before providing the details of this estimate, however, 
it may be helpful to sketch an outline of Bray's proof.

Given a real number $V \geq 0$, define
\begin{align}
\label{av}
A(V) = \inf_{\Omega}\{area(\partial \Omega) | vol(\Omega) = V \},
\end{align}
where $\Omega$ is any region in $M^3$, $vol(\Omega)$ is the volume of $\Omega$, and $area(\partial \Omega)$
is the $2$-dimensional surface area of the boundary.  Since $M^3$ is compact, there always exists a smooth
region whose boundary 
$\Sigma(V)$ attains the infimum $A(V)$.  Of course, $\Sigma(V)$ will necessarily 
have constant mean curvature.   

For a fixed value $V = V_0$ we consider a normal 
variation of $\Sigma(V_0)$, parametrized by the volume $V$. Let $A_{V_0}(V)$ denote the
area of this variation, and
primes denote differentiation with
respect to $V$.  Then $A_{V_0}'(V) = H$, where $H$ is the mean curvature of $\Sigma(V_0)$, and 
\begin{align}
\label{a21}
A_{V_0}(V_0)^2 A_{V_0}''(V_0) =  \int_{\Sigma(V_0)}[ -\|\Pi\|^2 - Ric(\nu,\nu)],
\end{align}
where $\Pi$ is the second fundamental form of $\Sigma(V_0)$ and $\nu$ is a unit normal.  From 
inequality (\ref{brayRc}) and
\begin{align}
\label{equalcase1}
\|\Pi\|^2 \geq \frac{1}{2}H^2,
\end{align}
we conclude
\begin{align}
\label{a0ppRc}
A_{V_0}''(V_0) \leq -\frac{1}{A_{V_0}(V_0)}\left(\frac{1}{2}A_{V_0}'(V_0)^2 + \epsilon Ric_0\right).
\end{align}
Since $\Sigma_{V_0}(V)$ may not attain the infimum in (\ref{av}), $A(V) \leq A_{V_0}(V)$.
Thus, $A(V)$ satisfies
\begin{align}
\label{appRc}
A''(V) \leq -\frac{1}{A(V)}\left(\frac{1}{2}A'(V)^2 + \epsilon Ric_0\right).
\end{align}

By the Gauss equation,
\begin{align*}
Ric(\nu,\nu) = \frac{1}{2}R - K + \frac{1}{2}H^2 - \frac{1}{2}\|\Pi\|^2,
\end{align*}
where $K$ is the Gauss curvature of $\Sigma(V_0)$.  Substituting this into 
(\ref{a21}) gives 
\begin{align}
\label{a22}
A_{V_0}(V_0)^2 A_{V_0}''(V_0) =  \int_{\Sigma(V_0)} 
[ -\frac{1}{2}R + K - \frac{1}{2}H^2 + \frac{1}{2}\|\Pi\|^2].
\end{align}
As Bray points out, the postivity of the Ricci curvature implies that $\Sigma(V_0)$ is connected.  Thus,
applying the Gauss-Bonnet formula and appealing to inequalities (\ref{Rlow0}) and
(\ref{equalcase1}) we get
\begin{align}
\label{a0ppR} 
A_{V_0}''(V_0) \leq \frac{4\pi}{A_{V_0}(V_0)^2} - \frac{1}{A_{V_0}(V_0)}\left( \frac{3}{4}A_{V_0}'(V_0)^2 + \frac{1}{2}R_0\right).
\end{align}
As before, since $A(V) \leq A_{V_0}(V)$ we have
\begin{align}
\label{appR}
A''(V) \leq \frac{4\pi}{A(V)^2} - \frac{1}{A(V)}\left( \frac{3}{4}A'(V)^2 + \frac{1}{2}R_0\right).
\end{align}

Next, let
\begin{align}
F(V) = A(V)^{3/2}.
\end{align}
By (\ref{appRc}) and (\ref{appR}), $F$ satisfies
\begin{align}
\label{FppRc}
F''(V) \leq -\frac{3\epsilon}{2}Ric_0F(V)^{-\frac{1}{3}},
\end{align}
\begin{align}
\label{FppR}
F''(V) \leq \frac{36\pi - F'(V)^2}{6F(V)} - \frac{3}{4}R_0F(V)^{-\frac{1}{3}}.  
\end{align}
Of course, one needs to properly interpret the sense in which these inequalities hold;
see (\cite{Bray}) for precise notions. 



Combining (\ref{FppRc}) and (\ref{FppR}), we have
\begin{align}
\label{mainode}
F''(V) \leq -\frac{1}{2}F^{-\frac{1}{3}}\max \{-\frac{36\pi - F'(V)^2}{3F(V)^{\frac{2}{3}}} + \frac{3}{2}R_0,3\epsilon Ric_0\}.
\end{align}

Consider the phase space associated to this differential inequality, which 
we view as the $xy$-plane with $x=F(V)$ and $y = F'(V)$.  Let $\gamma$ be
a path in phase space with intial value $V  = 0$ and terminal value
$V = \frac{1}{2}vol(M^3,g)$.  Then $\gamma$ starts at a point on the (positive)
$y$-axis and ends at a point on the (positive) $x$-axis.  
By (\ref{mainode}) this path must satisfy
the differential inequality
\begin{align}
\label{xyode}
\frac{dy}{dx} \leq -\frac{1}{2}x^{-\frac{1}{3}}y^{-1}\max\{-\frac{(36\pi-y^2)}{3x^{\frac{2}{3}}}+\frac{3}{2}R_0,3\epsilon Ric_0\}.
\end{align}
Also,
\begin{align}
\label{lineint}
\frac{1}{2}vol(M^3,g) = \int_{\gamma} dV = \int_{\gamma} \frac{dx}{y}.
\end{align}
A path which maximizes the line integral (\ref{lineint}) will be a path
which attains equality in (\ref{xyode}).  This results in an ODE which
can be explicitly solved, and by evaluating the line integral for this
path one obtains an upper estimate on the volume as in (\ref{volest}).

With this brief overview in mind, we now give an estimate of $\epsilon_0$. 

\begin{theorem}
\label{bestepsilon}
The constant $\epsilon_0 \leq \frac{1}{2}$.  
\end{theorem}

\begin{proof}
According to Bray's theorem, it suffices to show that $\alpha(\frac{1}{2}) = 1$; i.e., that
\begin{align}
\label{supless1}
\sup_{2\pi \leq z \leq 4\pi} \frac{1}{\pi^2}\left(
\begin{array}{rr} 
\int_0^{y(z)}(36\pi-\frac{27}{2}y(z)^{\frac{2}{3}}-\frac{9}{2}x^{\frac{2}{3}})^{-\frac{1}{2}}dx \\
+ \int_{y(z)}^{z^{\frac{3}{2}}}(36\pi-9y(z)x^{-\frac{1}{3}}-9x^{\frac{2}{3}})^{-\frac{1}{2}}dx\end{array}\right) = 1,
\end{align}
where
\begin{align}
\label{defy}
y = y(z) = z^{\frac{1}{2}}(4\pi-z).
\end{align}
To this end, let
\begin{align}
\label{defI1}
I_1(z) = \frac{1}{\pi^2} \int_0^y (36\pi-\frac{27}{2}y^{\frac{2}{3}}-\frac{9}{2}x^{\frac{2}{3}})^{-\frac{1}{2}}dx,
\end{align}
\begin{align}
\label{defI2}
I_2(z) = \frac{1}{\pi^2} \int_{y}^{z^{\frac{3}{2}}}(36\pi-9yx^{-\frac{1}{3}}-9x^{\frac{2}{3}})^{-\frac{1}{2}}dx.
\end{align}
We want to show that for $z \in [2\pi,4\pi]$, 
\begin{align}
\label{sumless1}
I_1(z) + I_2(z) \leq 1.
\end{align}
The first integral in (\ref{sumless1}) can be evaluated in closed form.  The second integral can be expressed in terms
of elliptic functions, though the resulting formula seems difficult to estimate.  Instead of this approach, we
will perform a change of variable and approximate the new integrand by one which can also be evaluated in closed
form.  It turns out to be much easier estimating both integrals in terms of this new variable; for this reason
we begin by analyzing $I_2$, where the substitution originates.  

Let $x=t^3$; then 
\begin{align}
\begin{split}
I_2 &= \frac{1}{\pi^2} \int_{y^{\frac{1}{3}}}^{z^{\frac{1}{2}}} [36\pi - 9yt^{-1}-9t^2]^{-\frac{1}{2}}(3t^2)dt \\
&= \frac{1}{\pi^2} \int_{y^{\frac{1}{3}}}^{z^{\frac{1}{2}}} \frac{t^{\frac{5}{2}}dt}{\sqrt{4\pi t - y - t^3}}.
\end{split}
\end{align}
Note that the denominator factors:
\begin{align*}
4 \pi t - y - t^3 = (z^{\frac{1}{2}}-t)(t^2 + z^{\frac{1}{2}}t - (4\pi-z)).
\end{align*}
Therefore, 
\begin{align*}
I_2 = \frac{1}{\pi^2} \int_{y^{\frac{1}{3}}}^{z^{\frac{1}{2}}} \frac{t^{\frac{5}{2}}dt}
{\sqrt{(z^{\frac{1}{2}}-t)(t^2 + z^{\frac{1}{2}}t - (4\pi-z))}}.
\end{align*}
Now perform another change of variable: let $s = tz^{-\frac{1}{2}}$; then
\begin{align*}
I_2 = \frac{z}{\pi^2} \int_{y^{\frac{1}{3}}z^{-\frac{1}{2}}}^{1} \frac{s^{\frac{5}{2}}ds}{\sqrt{(1-s)(s^2+s-(\frac{4\pi-z}{z}))}}.
\end{align*}
Let $\p = y^{\frac{1}{3}}z^{-\frac{1}{2}}.$  By (\ref{defy}),
\begin{align}
\label{phiz}
\p^3 = (\frac{4\pi-z}{z}),
\end{align}
\begin{align}
\label{zphi}
z = \frac{4\pi}{1+\p^3}.
\end{align}
Therefore, 
\begin{align}
\label{I2phi}
I_2 = (\frac{4}{\pi})(\frac{1}{1+\p^3}) \int_{\p}^{1} \frac{s^{\frac{5}{2}}ds}{\sqrt{(1-s)(s^2+s-\p^3)}}.
\end{align}
Since $z$ is a decreasing function of $\p$, we can change variables and view
$I_1$ and $I_2$ as functions of $\p$ (instead of $z$). 
 Note that $2\pi \leq z \leq 4\pi$,
while $0 \leq \p \leq 1$.

By doing a simple substitution the first integral can be evaluated in closed form:
\begin{align}
\label{I1eval}
\begin{split}
I_1(z) &= \frac{1}{\pi^2} \int_0^y (36\pi-\frac{27}{2}y^{\frac{2}{3}}-\frac{9}{2}x^{\frac{2}{3}})^{-\frac{1}{2}}dx \\
&= \frac{\sqrt{2}}{2\pi^2}(8\pi - 3y^{\frac{2}{3}})\left[ \arcsin\left(\frac{y^{\frac{1}{3}}}{\sqrt{8\pi-3 y^{\frac{2}{3}}}}\right) - \frac{2 y^{\frac{1}{3}}
\sqrt{2\pi - y^{\frac{2}{3}}}}{ 8\pi - 3 y^{\frac{2}{3}}} \right].
\end{split}
\end{align}
In order to rewrite (\ref{I1eval}) in terms of $\p$, we neeed to first express $y$ in terms of $\p$.  By (\ref{defy}) and (\ref{zphi}),
\begin{align}
\label{yphi}
y = \frac{ (4\pi)^{\frac{3}{2}}\p^3 }{ (1+\p^3)^{\frac{3}{2}} }.
\end{align}
Substituting this into (\ref{I1eval}) and carrying out the obvious simplifications, the result is
\begin{align}
\label{I1phi}
\begin{split}
I_1(\p) = (\frac{4}{\pi})(\frac{1}{1+\p^3})(\frac{\sqrt{2}}{2}) \Bigg[&
 (2 + 2\p^3 - 3\p^2) \arcsin \left( \frac{\p} { (2 + 2\p^3 - 3\p^2)^{\frac{1}{2}} }\right) \\
& - (2 + 2\p^3 - 4\p^2)^{\frac{1}{2}} \p \Bigg]. 
\end{split}
\end{align}

To establish the inequality
\begin{align}
\label{goal}
I_1(\p) + I_2(\p) \leq 1 \mbox{ for } \p \in [0,1] 
\end{align}
we divide the interval $[0,1]$ into two parts.  This division, or something like it, seems
necessary, since the contribution of the two integrals in the sum above is different for
$\p$ near $0$ and $\p$ near $1$.  More precisely, $I_1(\p) \to 0$  and $I_2(\p) \to 1$ as $\p \to 0$,
while $I_1(\p) \to 1/\sqrt{2}$ and $I_2(\p) \to 0$ as $\p \to 1$.
Therefore, in the subsections which follow we derive our estimates first on the interval $[0,\frac{4}{5}]$, then
on $[\frac{4}{5},1]$.  

\subsubsection{Estimate from $0$ to $\frac{4}{5}$}

We begin with an estimate of $I_1$:

\begin{proposition}
\label{I_1on.8}
For $\p \in [0,\frac{4}{5}]$, 
$$
I_1 \leq (\frac{4}{\pi})(\frac{1}{1+\p^3})(\frac{61}{100}\p^3).
$$
\end{proposition}

\begin{proof}

The proof relies on a sharp estimate of the arcsin term in (\ref{I1phi}).

\begin{lemma}
If $\beta \in (0,1]$, then for $x \in [0,\beta]$
\begin{align}
\label{calcineq}
\arcsin x \leq x + mx^3,
\end{align}
where 
$$
m = \left( \frac{ \arcsin \beta - \beta}{\beta^3} \right).
$$
\end{lemma}

\begin{proof}
This is equivalent to the inequality 
$$
\theta \leq \sin \theta + m\sin^3 \theta, \quad \theta \in [0,\arcsin \beta].
$$
Let $f(\theta) = \sin \theta + m\sin^3 \theta - \theta$.  We want to see that
$f(\theta) \geq 0$ for $\theta \in [0,\arcsin \beta]$.  Note that $f(0) = 0$, $f(\arcsin \beta) = 0$.
Thus, to show that $f(\theta) \geq 0$ it suffices to show that $(i) f(\theta) > 0$ for $\theta > 0$ small, 
and $(ii) f^{\prime}$ has exactly one zero in the open interval $(0,\arcsin \beta)$.  Of course,
since $f(0)=f(\arcsin \beta)=0$
Rolle's theorem guarantees that $f^{\prime}(\theta_0) = 0$ for some $\theta_0 \in   (0,\arcsin \beta)$.

If we write out the Taylor expansion of $f$ near zero,  
$$
f(\theta) = (m - \frac{1}{6})\theta^3 + O(\theta^5).
$$
Thus, if we can show that $m > \frac{1}{6}$ then $(i)$ will follow.  To this end, define another function
$h(\beta) = \arcsin \beta - \beta - \frac{1}{6}\beta^3$.  Then 
$$
h^{\prime}(\beta) = \frac{1}{\sqrt{1-\beta^2}} - 1 - \frac{1}{2}\beta^2.
$$
It is easy to see that $h^{\prime}(\beta) > 0$ for $\beta \in (0,1)$: just differentiate again and use the fact that $h^{\prime}(0)=0$. 
Thus, $h(\beta) > h(0) = 0$ for $\beta \in (0,1)$, which implies $m > \frac{1}{6}$.

To prove $(ii)$, note 
$f^{\prime}(\theta) = (1+3m) \cos \theta - 3m \cos^3 \theta - 1$.
Let 
$$
p(z) = (1+3m)z - 3mz^3  - 1.
$$
If we can show that $p$ has exactly one zero in the interval $(\cos (\arcsin \beta),1) = (\sqrt{1-\beta^2},1)$, then $(ii)$ will follow.
Let $z_0 = \cos \theta_0$; then $p(z_0) =0$.  Also, $p(1) = 0$.  Thus, $p$ has two
zeros in the closed interval $[\sqrt{1-\beta^2},1]$: $z_0 \in (\sqrt{1-\beta^2},1)$, and $z_1 = 1$.
Since $p$ is a cubic polynomial, it must have a 
third zero $z_2$.  But notice 
$$
\lim_{z \to -\infty} p(z) = +\infty
$$
while $p(0) = -1$.  Consequently, $z_2 < 0$, and $p$ has only one zero in the
open interval $(\sqrt{1-\beta^2},1)$.  
\end{proof}

Using the preceding Lemma, we estimate the $\arcsin$ term in (\ref{I1phi}) as follows.  First, 
observe that 
\begin{align}
\label{biggerthan1}
2 + 2\p^3 - 3\p^2 \geq 1.
\end{align}
This follows from the fact that $\p^2 \leq \frac{2}{3}\p^3 + \frac{1}{3}$, and hence 
$-3\p^2 \geq -2\p^3 -1$.  A consequence of (\ref{biggerthan1}) is that
\begin{align*}
\frac{\p}{ (2 + 2\p^3 - 3\p^2)^{\frac{1}{2}} } \leq \p.
\end{align*}
Therefore, $x \equiv \p/(2 + 2\p^3 - 3\p^2)^{\frac{1}{2}} \in [0,\frac{4}{5}]$
whenever $\p \in [0,\frac{4}{5}]$.  From (\ref{calcineq}) we conclude 
\begin{align*}
\arcsin \left( \frac{\p} { (2 + 2\p^3 - 3\p^2)^{\frac{1}{2}} }\right) 
\leq \frac{\p}{ (2 + 2\p^3 - 3\p^2)^{\frac{1}{2}} } 
+ m_0 \frac{\p^3} { (2 + 2\p^3 - 3\p^2)^{\frac{3}{2}} },
\end{align*}
where
\begin{align}
\label{m0}
m_0 = \left( \frac{\arcsin \frac{4}{5} - \frac{4}{5} }{ (\frac{4}{5})^3 }\right).
\end{align}
Therefore,
\begin{align*}
\begin{split}
(2 + 2\p^3 - 3\p^2) &\arcsin \left( \frac{\p} { (2 + 2\p^3 - 3\p^2)^{\frac{1}{2}} }\right) 
- (2 + 2\p^3 - 4\p^2)^{\frac{1}{2}}\p \\
&\leq (2 + 2\p^3 - 3\p^2)^{\frac{1}{2}}\p  + 
m_0 \frac{\p^3} { (2 + 2\p^3 - 3\p^2)^{\frac{1}{2}} } 
 - (2 + 2\p^3 - 4\p^2)^{\frac{1}{2}} \p \\
&=  \left[ (2 + 2\p^3 - 3\p^2)^{\frac{1}{2}}  - 
(2 + 2\p^3 - 4\p^2)^{\frac{1}{2}} \right] \p 
+ m_0 \frac{\p^3} { (2 + 2\p^3 - 3\p^2)^{\frac{1}{2}} } 
\end{split}
\end{align*}
whenever $\p \in [0,\frac{4}{5}]$. 
For the first term above,
\begin{align*}
\begin{split}
& (2 + 2\p^3 - 3\p^2)^{\frac{1}{2}} - (2 + 2\p^3 - 4\p^2)^{\frac{1}{2}} \\
&= \frac{ \left[ (2 + 2\p^3 - 3\p^2)^{\frac{1}{2}} - (2 + 2\p^3 - 4\p^2)^{\frac{1}{2}}\right]
\left[ (2 + 2\p^3 - 3\p^2)^{\frac{1}{2}} + (2 + 2\p^3 - 4\p^2)^{\frac{1}{2}}\right] }
{ \left[ (2 + 2\p^3 - 3\p^2)^{\frac{1}{2}} + (2 + 2\p^3 - 4\p^2)^{\frac{1}{2}}\right] } \\
&= \frac { \p^2 } { \left[ (2 + 2\p^3 - 3\p^2)^{\frac{1}{2}} + (2 + 2\p^3 - 4\p^2)^{\frac{1}{2}}\right] }.
\end{split}
\end{align*}
On the interval $[0,\frac{4}{5}]$, $l(\p) = \frac{ (2 + 2\p^3 - 4\p^2) }{ (2 + 2\p^3 - 3\p^2) }$ is
decreasing; thus $ (2 + 2\p^3 - 4\p^2)^{\frac{1}{2}} \geq l(\frac{4}{5})^{\frac{1}{2}} (2 + 2\p^3 - 3\p^2)^{\frac{1}{2}}$.
Substituting this into the expression above, 
\begin{align*}
(2 + 2\p^3 - 3\p^2)^{\frac{1}{2}} - (2 + 2\p^3 - 4\p^2)^{\frac{1}{2}} \leq \frac{ \p^2 }
 { (1 + l(\frac{4}{5})^{\frac{1}{2}} )(2 + 2\p^3 - 3\p^2)^{\frac{1}{2}} }.
\end{align*}
Therefore,
\begin{align*}
\begin{split}
(2 + 2\p^3 - 3\p^2) \arcsin & \left( \frac{\p} { (2 + 2\p^3 - 3\p^2)^{\frac{1}{2}} }\right) 
- (2 + 2\p^3 - 4\p^2)^{\frac{1}{2}} \p \\ 
& \leq \left( \frac{1}{ 1+ l(\frac{4}{5})^{\frac{1}{2} }} + m_0\right )\frac{ \p^3 } { (2 + 2\p^3 - 3\p^2)^{\frac{1}{2}} }. 
\end{split}
\end{align*}
Substituting this into (\ref{I1phi}) we conclude
\begin{align*}
\begin{split}
I_1(\p) & \leq (\frac{4}{\pi})(\frac{1}{1+\p^3})(\frac{\sqrt{2}}{2})\left(\frac{1}{ 1+ l(\frac{4}{5})^{\frac{1}{2}} }  + m_0 \right)\p^3 \\
& = (\frac{4}{\pi})(\frac{1}{1+\p^3})c_0 \p^3,
\end{split}
\end{align*}
where
\begin{align*}
\begin{split}
c_0 & = (\frac{\sqrt{2}}{2})\left( \frac{1}{ 1+ l(\frac{4}{5})^{\frac{1}{2}} } + m_0 \right) \\
&= (\frac{\sqrt{2}}{2})\Big[ \frac{1}{ 1+ (\frac{29}{69})^{\frac{1}{2}} } +  
\left( \frac{\arcsin \frac{4}{5} - \frac{4}{5} }{ (\frac{4}{5})^3 }\right) \Big] \\
&= 0.604795... \\
&< \frac{61}{100}.
\end{split}
\end{align*}
\end{proof}

To estimate $I_2$, we begin by rewriting the integrand in (\ref{I2phi}):
\begin{align}
\label{deff}
\begin{split}
\frac{s^{\frac{5}{2}}}{\sqrt{(1-s)(s^2+s-\p^3)}} &= \frac{s^{\frac{5}{2}}}{\sqrt{(1-s)(s^2+s)}}\sqrt{\frac{s^2+s}{s^2+s-\p^3}} \\
&= \frac{s^2}{\sqrt{(1-s)(1+s)}}f(s),
\end{split}
\end{align}
where
\begin{align*}
f(s) = \frac{(s^2+s)^{\frac{1}{2}}}{(s^2+s-\p^3)^{\frac{1}{2}}}.
\end{align*}
Differentiating,
$$
f'(s)=-\frac{1}{2}\p^3(2s+1)(s^2 +s-\p^3)^{-\frac{3}{2}}(s^2+s)^{-\frac{1}{2}}.
$$
Since
$$
s^2+s-\p^3\leq s^2 + s,
$$
it follows 
$$
(s^2+s-\p^3)^{-\frac{3}{2}} \geq (s^2+s)^{-\frac{3}{2}}.
$$
Therefore, 
$$
f'(s) \leq -\frac{1}{2}\p^3\frac{(2s+1)}{(s^2+s)^2}.
$$
By the fundamental theorem of calculus,
\begin{align*}
\begin{split}
f(s) - f(\p) &\leq \int_{\p}^{s} -\frac{1}{2}\p^3\frac{(2x+1)}{(x^2+x)^2} dx \\
&= \frac{1}{2}\p^3(x^2+x)^{-1}\vert_{x=\p}^{x=s} \\
&= \frac{1}{2}\p^3\left( \frac{1}{s^2+s} - \frac{1}{\p^2+\p}\right),
\end{split}
\end{align*}
hence
$$
f(s) \leq \left(f(\p) - \frac{1}{2}\frac{\p^2}{1+\p}\right) + \frac{1}{2}\p^3\frac{1}{s(1+s)}.
$$
Substituting this inequality into (\ref{deff}) we have
\begin{align}
\label{2ints}
\begin{split}
\int_{\p}^{1} \frac{s^2}{\sqrt{(1-s)(1+s)}} & f(s) ds  \\ 
&\leq \int_{\p}^{1} \frac{s^2}{\sqrt{(1-s)(1+s)}} 
\left[ \left(f(\p) - \frac{1}{2}\frac{\p^2}{1+\p}\right) + \frac{1}{2}\p^3\frac{1}{s(1+s)} \right] ds \\
&= \left(f(\p) - \frac{1}{2}\frac{\p^2}{1+\p}\right) \int_{\p}^{1} \frac{s^2}{\sqrt{1-s^2}}ds 
+ \frac{1}{2}\p^3 \int_{\p}^{1} \frac{s}{\sqrt{(1-s)(1+s)^3}} ds. 
\end{split}
\end{align}
Both integrals in (\ref{2ints}) are elementary:
\begin{align*}
\begin{split}
\left(f(\p) - \frac{1}{2}\frac{\p^2}{1+\p}\right) \int_{\p}^{1} \frac{s^2}{\sqrt{1-s^2}} & ds \\
&= \left(f(\p) - \frac{1}{2}\frac{\p^2}{1+\p}\right)\left[-\frac{1}{2}s\sqrt{1-s^2} + \frac{1}{2}\arcsin s \right]_{s=\p}^{s=1} \\
&= \left(f(\p) - \frac{1}{2}\frac{\p^2}{1+\p}\right) \left[ \frac{\pi}{4} + \frac{1}{2}\p\sqrt{1-\p^2}-\frac{1}{2}\arcsin \p \right],\\
\frac{1}{2}\p^3 \int_{\p}^{1} \frac{s}{\sqrt{(1-s)(1+s)^3}} & ds 
= \frac{1}{2}\p^3 \left[ \sqrt{\frac{1-s}{1+s}} + \arcsin s \right]_{s=\p}^{s=1} \\
&= \frac{1}{2}\p^3 \left[ \frac{\pi}{2} - \sqrt{\frac{1-\p}{1+\p}} - \arcsin \p \right].
\end{split}
\end{align*}
Combining the above and substituting into (\ref{I2phi}) we get
\begin{align}
\label{I2leq}
\begin{split}
I_2 \leq &  (\frac{4}{\pi})(\frac{1}{1+\p^3})\Big\{ \left(\frac{(1+\p)^{\frac{1}{2}}}{(1+\p-\p^2)^{\frac{1}{2}}} - \frac{1}{2}\frac{\p^2}{1+\p}\right)
\left[ \frac{\pi}{4} + \frac{1}{2}\p\sqrt{1-\p^2}-\frac{1}{2}\arcsin \p \right] \\ 
& + \p^3 \left[ \frac{\pi}{2} - \sqrt{\frac{1-\p}{1+\p}} - \arcsin \p \right] \Big\} \\
&= (\frac{4}{\pi})(\frac{1}{1+\p^3})\Big \{ E(\p)F(\p) + G(\p)\p^3 \Big\},
\end{split}
\end{align}
where
\begin{align}
\label{defE}
E(\p) = \frac{(1+\p)^{\frac{1}{2}}}{(1+\p-\p^2)^{\frac{1}{2}}} - \frac{1}{2}\frac{\p^2}{1+\p},
\end{align}
\begin{align}
\label{defF}
F(\p) = \frac{\pi}{4} + \frac{1}{2}\p\sqrt{1-\p^2}-\frac{1}{2}\arcsin \p,
\end{align}
\begin{align}
\label{defJ}
G(\p) = \frac{\pi}{4} - \frac{1}{2}\sqrt{\frac{1-\p}{1+\p}} - \frac{1}{2}\arcsin \p.
\end{align}

\begin{proposition}
\label{I2intermsofp3}
For $\p \in [0,\frac{4}{5}]$,
\begin{align}
\label{I2on0to.8}
I_2 \leq (\frac{4}{\pi})(\frac{1}{1+\p^3})\Big\{ \frac{\pi}{4} + (\frac{\pi}{16} - \frac{1}{3})\p^3 
 + (\frac{\pi}{4} - \frac{1}{2})\p^3  \Big\}.
\end{align}
\end{proposition}

\begin{proof}
The proof of (\ref{I2on0to.8}) is based on the following estimates of 
$E$,$F$, and $G$.

\begin{lemma}
\label{2.2}

\noindent
$(i)$ For $\p \in [0,\frac{2}{5}]$, 
\begin{align}
\label{2.2resulti}
E(\p) \leq 1 + \frac{1}{2}\p^4.
\end{align}

\noindent
$(ii)$  For $\p \in [\frac{2}{5},\frac{4}{5}]$, 
\begin{align}
\label{2.2resultii}
E(\p) \leq 1 + \frac{125}{434}\p^4.
\end{align}
\end{lemma}

\begin{proof}
First, write
$$
E(\p) = E_1(\p) + E_2(\p),
$$ 
where
\begin{align*}
E_1(\p) = \frac{(1+\p)^{\frac{1}{2}}}{(1+\p-\p^2)^{\frac{1}{2}}},
\end{align*}
\begin{align*}
E_2(\p) = - \frac{1}{2}\frac{\p^2}{1+\p}.
\end{align*}
Then
\begin{align*}
\begin{split}
E_1(\p) &= \frac{(1+\p)^{\frac{1}{2}}}{(1+\p-\p^2)^{\frac{1}{2}}} - 1 + 1 \\
&= \frac{(1+\p)^{\frac{1}{2}} - (1+\p - \p^2)^{\frac{1}{2}}}{(1+\p-\p^2)^{\frac{1}{2}}} + 1 \\
&= \frac{\left[ (1+\p)^{\frac{1}{2}} - (1+\p - \p^2)^{\frac{1}{2}}\right]  \left[ (1+\p)^{\frac{1}{2}} + (1+\p - \p^2)^{\frac{1}{2}}\right]}
{ (1+\p - \p^2)^{\frac{1}{2}} \left[  (1+\p)^{\frac{1}{2}} + (1+\p - \p^2)^{\frac{1}{2}}\right]} + 1 \\
&= \frac{ \p^2 }{ (1+\p - \p^2)^{\frac{1}{2}} \left[  (1+\p)^{\frac{1}{2}} + (1+\p - \p^2)^{\frac{1}{2}}\right]} + 1.
\end{split}
\end{align*}
Since $(1+\p)^{\frac{1}{2}} \geq (1+\p-\p^2)^{\frac{1}{2}}$, we can estimate the denominator above by 
$$
(1+\p - \p^2)^{\frac{1}{2}} \left[ (1+\p)^{\frac{1}{2}} + (1+\p - \p^2)^{\frac{1}{2}}\right] \geq 2(1+\p-\p^2).
$$
Thus,
$$
E_1(\p) \leq \frac{\p^2}{2(1+\p-\p^2)} + 1.
$$
So
\begin{align}
\label{casesiii}
\begin{split}
E(\p) &= E_1(\p) + E_2(\p) \\
&\leq \frac{\p^2}{2(1+\p-\p^2)} + 1 - \frac{\p^2}{2(1+\p)} \\
&= 1 + \frac{\p^4}{2(1+\p)(1+\p-\p^2)} \\
&\equiv 1 + \frac{\p^4}{D(\p)}, 
\end{split}
\end{align}
where 
\begin{align*}
D(\p) = 2(1+\p)(1+\p-\p^2).
\end{align*}
Differentiating, we see that $D^{\prime}(\p) = 4 - 6\p^2 > 0$ for $\p \in [0,\frac{4}{5}]$.   Thus,
on the interval $[0, \frac{2}{5}]$ we have $D(\p) \geq D(0) = 2$, while on the interval $[\frac{2}{5},\frac{4}{5}]$
we have $D(\p) \geq D(\frac{2}{5}) = \frac{434}{125}$.  Substituting these inequalities into (\ref{casesiii}) we obtain
(\ref{2.2resulti}) and (\ref{2.2resultii}).  
\end{proof}

\begin{lemma}
\label{2.4}
For $\p \in [0,\frac{4}{5}]$,
\begin{align}
\label{2.4result}
F(\p) \leq \frac{\pi}{4} - \frac{1}{3}\p^3.
\end{align}
\end{lemma}

\begin{proof}
Since $F(0) = \frac{\pi}{4}$ and $F^{\prime}(\p) = \frac{-\p^2}{\sqrt{1-\p^2}} \leq -\p^2$, upon integrating we find
\begin{align*}
\begin{split}
F(\p) - F(0) &= \int_0^{\p} F^{\prime}(s) ds  \\
&\leq \int_0^{\p} -s^2 ds \\
&= -\frac{1}{3}\p^3. 
\end{split}
\end{align*}
\end{proof}

\begin{lemma}
\label{2.6}
For $\p \in [0,\frac{4}{5}]$, 
\begin{align}
\label{2.6result}
E(\p)F(\p) \leq \frac{\pi}{4} + (\frac{\pi}{16} - \frac{1}{3})\p^3.
\end{align}
\end{lemma}

\begin{proof}
When $\p \in [0,\frac{2}{5}]$, by (\ref{2.2resulti}) and (\ref{2.4result})
\begin{align*}
\begin{split}
E(\p)F(\p) & \leq \left( 1 + \frac{1}{2}\p^4 \right) \left( \frac{\pi}{4} - \frac{1}{3}\p^3 \right) \\
& = \frac{\pi}{4} - \frac{1}{3}\p^3 + \frac{\pi}{8}\p^4 - \frac{1}{6}\p^7 \\
& \leq  \frac{\pi}{4} - \frac{1}{3}\p^3 + \frac{\pi}{8}\p^4.  
\end{split}
\end{align*}
Since $\p \leq \frac{2}{5}$ implies $\frac{\pi}{8}\p^4 \leq \frac{\pi}{20}\p^3 \leq \frac{\pi}{16}\p^3$, 
(\ref{2.6result}) follows.

Similarly, for $\p \in [\frac{2}{5},\frac{4}{5}]$, by (\ref{2.2resultii}) and (\ref{2.4result})
\begin{align*}
\begin{split}
E(\p)F(\p) & \leq \left( 1 + \frac{125}{434}\p^4 \right) \left( \frac{\pi}{4} - \frac{1}{3}\p^3 \right) \\
& = \frac{\pi}{4} - \frac{1}{3}\p^3 + \frac{125\pi}{(4)(434)}\p^4 - \frac{125}{(3)(434)}\p^7 \\
& \leq  \frac{\pi}{4} - \frac{1}{3}\p^3 + \frac{125\pi}{(4)(434)}\p^4.  
\end{split}
\end{align*}
When $\p \leq \frac{4}{5}$,  $\frac{125\pi}{(4)(434)}\p^4 \leq \frac{25\pi}{434}\p^3 \leq \frac{\pi}{16}\p^3$, and once
again (\ref{2.6result}) holds.  
\end{proof}

\begin{lemma}
\label{2.5}
For $\p \in [0,\frac{4}{5}]$,
\begin{align}
\label{2.5result}
G(\p) \leq \frac{\pi}{4} - \frac{1}{2}.
\end{align}
\end{lemma}

\begin{proof}
Since
$$
G^{\prime}(\p) = -\frac{\p}{\sqrt{(1-\p)(1+\p)}} \leq 0,
$$
it follows that $G(\p) \leq G(0) = \frac{\pi}{4} - \frac{1}{2}$.
\end{proof}

To complete the proof of Proposition \ref{I2intermsofp3}, notice that (\ref{I2on0to.8}) follows from 
(\ref{I2leq}), (\ref{2.6result}), and (\ref{2.5result}).  
\end{proof}

Combining the results of Propositions \ref{I_1on.8} and \ref{I2intermsofp3}, we see that  
\begin{align*}
I_1 + I_2  \leq (\frac{4}{\pi})(\frac{1}{1+\p^3})\Big\{ \frac{\pi}{4} + (\frac{\pi}{16} - \frac{1}{3})\p^3 
 + (\frac{\pi}{4} - \frac{1}{2})\p^3 + \frac{61}{100}\p^3 \Big\}.
\end{align*}
Therefore,
$$
I_1 + I_2 - 1 \leq (\frac{4}{\pi})(\frac{\p^3}{1+\p^3}) \left[\frac{\pi}{16} - \frac{5}{6} + \frac{61}{100} \right] \leq 0.
$$
It follows that $I_1 + I_2 \leq 1$ for $\p \in [0,\frac{4}{5}]$.

\subsubsection{Estimate from $\frac{4}{5}$ to $1$}

When $\p \in [\frac{4}{5},1]$, we need to use different estimates of $I_1$ and $I_2$.
First, recall formula (\ref{I1phi}): 
\begin{align*}
\begin{split}
I_1(\p) = (\frac{4}{\pi})(\frac{1}{1+\p^3})(\frac{\sqrt{2}}{2}) \Bigg[&
 (2 + 2\p^3 - 3\p^2) \arcsin \left( \frac{\p} { (2 + 2\p^3 - 3\p^2)^{\frac{1}{2}} }\right) \\
& - (2 + 2\p^3 - 4\p^2)^{\frac{1}{2}} \p \Bigg]. 
\end{split}
\end{align*}
By (\ref{biggerthan1})
\begin{align*}
2 - 3 \p^2 + 2 \p^3 \geq 1,
\end{align*} 
and since $\arcsin$ is increasing, it follows that
\begin{align}
\label{I1est}
I_1(\p) \leq (\frac{4}{\pi})(\frac{1}{1+\p^3})(\frac{\sqrt{2}}{2}) \Bigg[
 (2 + 2\p^3 - 3\p^2) \arcsin \p  - (2 + 2\p^3 - 4\p^2)^{\frac{1}{2}} \p \Bigg]. 
\end{align}

To estimate $I_2$ we use the fact that $\p \leq s$ in the integrand in $I_2$, so
\begin{align*}
s^2 + s - \p^3 \geq s^2 + s -s \p^2 = s(s + 1 - \p^2) \geq s^2.
\end{align*}
Therefore,
\begin{align*}
\int_{\p}^1 \frac{s^{5/2} ds}
{ \sqrt{ (1-s)(s^2 + s - \p^3)}}
& \leq \int_{\p}^1 \frac{s^{3/2} ds}
{ \sqrt{ (1-s)}}
\leq  \int_{\p}^1 \frac{s ds}
{ \sqrt{ (1-s)}}\\
&= {\sqrt{1 - s}}\,\left( - \frac{4}{3}   - 
\frac{2}{3}s \right) \Bigg|_\p^1
=  {\sqrt{1 - \p}}\,\left(  \frac{4}{3}   +
\frac{2}{3}\p \right). 
\end{align*}
Substituting this into (\ref{I2phi}) we conclude 
\begin{align}
\label{I2est}
I_2 \leq 
 \frac{4}{\pi} \frac{1}{1 + \p^3} {\sqrt{1 - \p}}\,\left(  \frac{4}{3}   +
\frac{2}{3}\p \right). 
\end{align}

Combining (\ref{I1est}) and (\ref{I2est}), and observing 
that $1 -2 \p^2 + \p^3 = (1-\p)(1 + \p - \p^2)$, we obtain
\begin{align}
\label{thisisH}
\begin{split}
I_1(\p) + I_2(\p) 
&\leq H(\p) \equiv (\frac{4}{\pi})(\frac{1}{( 1 + \p^3)}  \Bigg\{
\frac{\sqrt{2}}{2} \left( 2 - 3 \p^2 + 2 \p ^3 \right)
     \arcsin \p \\
& + \sqrt{ 1 - \p} 
\left( - \p \sqrt{ 1 + \p - \p^2}
+ \frac{4}{3} +\frac{2}{3}\p \right) \Bigg\}. 
\end{split}
\end{align}
It is elementary to estimate 
that $H(\frac{4}{5}) < .9881 < 1$. We will show that 
for $\p \in [\frac{4}{5},1]$, $H'( \p) < 0$, and
therefore $I_1(\p) + I_2(\p) 
\leq H(\p) < 1$.

 A computation shows that 
\begin{align}
\label{H'}
\begin{split}
H'(\p) & =  \frac{1}{\pi}\Bigg\{
 \frac{2 {\sqrt{2}} \left( -6 \p + 6 \p^2 \right) 
\arcsin \p}{1 + \p^3} - \frac{6 {\sqrt{2}}\p^2
 \left( 2 - 3\p^2 + 2\p^3 \right) \arcsin \p}
  {{\left( 1 + \p^3 \right) }^2}\\
& - \frac{2\left(  \frac{4}{3} + \frac{2\p}{3} \right) }
   {{\sqrt{1 - \p}}\left( 1 + \p^3 \right)  }    
 - \frac{ 2\p\left( -4\p + 3\p^2 \right) }
      {\left( 1 + \p^3 \right) { \sqrt{1 - \p} 
\sqrt{1 + \p - \p^2}}} 
 + \frac{2{\sqrt{2}}\left( 2 - 3\p^2 + 2\p^3 \right) }
      {{\sqrt{1 - \p^2}}\left( 1 + \p^3 \right) } \\
& \frac{-12{\sqrt{1 - \p}}\left( \frac{4}{3} + \frac{2\p}{3} \right) \p^2}
   {{\left( 1 + \p^3 \right) }^2  } 
+  \frac{12\p^3{\sqrt{1-\p} \sqrt{1 + \p - \p^2}}}{
{\left( 1 + \p^3 \right) }^2}\\
& + \frac{8{\sqrt{1 - \p}}}{3\left( 1 + \p^3 \right)  } -
 \frac{4{ \sqrt{1 -\p}\sqrt{1 +\p - \p^2}}}{1 + \p^3} \Bigg\}.
\end{split}
\end{align}
To see that $H^{\prime} < 0$ we will estimate each line of (\ref{H'}).
First we observe that the arcsin terms simplify to
\begin{align}
\frac{6 \sqrt{2} \p}{\pi (1 + \p^3)^2} (\p^3 -2) \arcsin{\p}.
\end{align}
\begin{lemma}
For $ \p \in [\frac{1}{2},1]$, 
\begin{align*}
& \frac{1}{\pi} \Bigg(- \frac{2\left( \frac{4}{3} + \frac{2\p}{3} 
\right) }
   {{\sqrt{1 - \p}}\left( 1 + \p^3 \right)  } -   
\frac{ 2\p\left( -4\p + 3\p^2 \right) }
      {\left( 1 + \p^3 \right) { \sqrt{1 - \p} 
\sqrt{1 + \p - \p^2}}} 
 + \frac{2{\sqrt{2}}\left( 2 - 3\p^2 + 2\p^3 \right) }
      {{\sqrt{1 - \p^2}}\left( 1 + \p^3 \right)} \Bigg) \\
& \leq 
\frac{2}{\pi} \frac{1}{(1 + \p^3)} 
( \frac{2}{3} + 2 \sqrt{2} ) \sqrt{ 1 - \p}.
\end{align*}
\end{lemma}
\begin{proof}
We begin by rewriting the left hand side  as 
\begin{align*}
\frac{2}{\pi} \frac{1}{ \sqrt{ 1 - \p}(1 + \p^3)} J(\p),
\end{align*}
where 
\begin{align*}
J(\p) = -\frac{4}{3}   - \frac{2}{3}\p + 
   \frac{\p\left( 4\p - 3\p^2 \right) }{{\sqrt{1 + \p - \p^2}}} + 
   \frac{{\sqrt{2}}\left( 2 - 3\p^2 + 2\p^3 \right) }{{\sqrt{1 + \p}}}.
\end{align*}
The polynomial $4\p^2 - 3\p^3 + \p - 2 \leq 0$ 
for $ \p \in [0,1]$; therefore,
\begin{align*}
 \frac{\p\left( 4\p - 3\p^2 \right) }{{\sqrt{1 + \p - \p^2}}}
\leq  4\p^2 - 3\p^3 \leq 2 - \p.
\end{align*}
Next we use the inequalities 
\begin{align}
\frac { \sqrt{2}}{ \sqrt{1+\p}}
\leq 1 + ( \sqrt{2} -1)(1-\p) , \ \p \in [0,1],
\end{align}
and
\begin{align}
2 - 3\p^2 + 2 \p^3 \leq 1 + (1-\p), \ \p \in [\frac{1}{2},1].
\end{align}
To derive these inequalities, simply use the fact that 
a convex function lies below the line segment between
the endpoints. 
It follows that
\begin{align*} 
 \frac{{\sqrt{2}}\left( 2 - 3\p^2 + 2\p^3 \right) }{{\sqrt{1 + \p}}}
&\leq ( 1 + (\sqrt{2}-1)(1-\p)) ( 1 + (1 -\p)) \\
&\leq 1 + \sqrt{2} (1 - \p) + (\sqrt{2}-1) (1 - \p)^2
\leq  1 + (2 \sqrt{2}-1) (1 - \p).
\end{align*}
Combining the preceding estimates, we obtain 
\begin{align*}
J( \p) \leq -\frac{4}{3}   - \frac{2}{3}\p + 2 - \p + 
1 + (2 \sqrt{2} -1)(1 - \p)
= ( \frac{2}{3} + 2 \sqrt{2})(1-\p). 
\end{align*}
\end{proof}
\begin{lemma} For $\p \in [\frac{4}{5},1]$,  
\begin{align*}
& \frac{1}{\pi} \Bigg( \frac{-12{\sqrt{1 - \p}}\left( \frac{4}{3} + \frac{2\p}{3} \right) \p^2}
   {{\left( 1 + \p^3 \right) }^2  } 
+  \frac{12\p^3{\sqrt{1-\p} \sqrt{1 + \p - \p^2}}}{
{\left( 1 + \p^3 \right) }^2} \Bigg) 
 \leq -11 \frac{ \p^2 \sqrt{ 1 - \p}}{\pi (1 + \p^3)^2}.
\end{align*}
\end{lemma}
\begin{proof}
Write the left hand side as 
\begin{align*}
\frac{ 12 \p^2 \sqrt{ 1 - \p}}{\pi (1 + \p^3)^2} K(\p),
\end{align*}
where 
\begin{align*}
K(\p) = - \frac{4}{3} - \frac{2}{3} \p + \p \sqrt{ 1 + \p - \p^2}.
\end{align*} 
It is easy to verify that $K$ is a concave function 
for $0 \leq \p \leq 1$, and therefore $K$ lies below 
its tangent line at $1$. A computation shows 
$K'(1) = -\frac{1}{6}$, so 
$K(\p) \leq -1 + \frac{1}{6} (1 -\p)$. 
Then $K(\p) \leq -\frac{29}{30}
< - \frac{11}{12}$ for $ \p \geq \frac{4}{5}$.
\end{proof}
\begin{lemma}
For $\p \in [0,1]$, 
\begin{align*}
\frac{1}{\pi} \Bigg(  \frac{8{\sqrt{1 - \p}}}{3\left( 1 + 
\p^3 \right)  } -
 \frac{4{ \sqrt{1 -\p}\sqrt{1 +\p - \p^2}}}{1 + \p^3} 
\Bigg) \leq - \frac{4}{3\pi} \frac{ \sqrt{ 1 - \p}}{1 + \p^3}. 
\end{align*}
\end{lemma}
\begin{proof}
Write the left hand side as 
\begin{align*}
\frac{4}{\pi} \frac{ \sqrt{ 1 - \p}}{1 + \p^3} \left(
\frac{2}{3} - \sqrt{ 1 + \p - \p^2} \right).
\end{align*}
The function $ \frac{2}{3} - \sqrt{ 1 + \p - \p^2} $
is clearly convex for $\p \in [0,1]$, therefore it achieves
its maximum at the endpoints, where it equals $-\frac{1}{3}$. 
\end{proof}
Combining the preceding Lemmas, we have the estimate
\begin{align*}
H'(\p) &\leq 
\frac{6 \sqrt{2} \p}{\pi (1 + \p^3)^2} (\p^3 -2) \arcsin{\p}
+ 
\frac{2}{\pi} \frac{1}{(1 + \p^3)} 
( \frac{2}{3} + 2 \sqrt{2} ) \sqrt{ 1 - \p}
 -11 \frac{ \p^2 \sqrt{ 1 - \p}}{\pi (1 + \p^3)^2}\\
& \ \ \ \ \ \ \ \ - \frac{4}{3\pi} \frac{ \sqrt{ 1 - \p}}{1 + \p^3} \\
&\leq \frac{1}{\pi (1 + \p^3)^2} 
\Bigg( 6 \sqrt{2} \p (\p^3 -2) \arcsin{\p}
+ \sqrt{ 1 - \p} \left( 
  4 \sqrt{2}(1 + \p^3) 
- 11 \p^2 \right) \Bigg).
\end{align*}
The polynomial $ 4 \sqrt{2}(1 + \p^3) 
- 11 \p^2$ is decreasing on $[0,1]$, so 
\begin{align*}
4 \sqrt{2}(1 + \p^3) 
- 11 \p^2 \leq 4 \sqrt{2}(1 + (\frac{4}{5})^3) 
- 11 (\frac{4}{5})^2 < 2  \mbox{ for } \p \in [\frac{4}{5},1].  
\end{align*}
Furthermore, $\arcsin\p > \frac{5}{6}$ for $ \p \in  [\frac{4}{5},1]$,
so we have 
\begin{align*}
H'(\p) &\leq \frac{1}{\pi (1 + \p^3)^2} 
\Bigg( 5 \sqrt{2} \p (\p^3 -2) 
+ 2  \Bigg).
\end{align*}
A simple calculation shows that the polynomial $ \p( \p^3-2)$ 
is increasing on $[\frac{4}{5},1]$, so 
$\p ( \p^3 -2) \leq -1$ for $\p \in [\frac{4}{5},1]$. 

Finally, by combining the above estimates we have 
\begin{align}
H'(\p) &\leq \frac{1}{\pi (1 + \p^3)^2} 
\Bigg( -5 \sqrt{2} + 2 \Bigg) < 0 \mbox{ for } \p
 \in [\frac{4}{5},1].
\end{align}
This completes the proof of Theorem \ref{bestepsilon}.

\end{proof}

\subsection{n = 4} 

In four dimensions our estimate of the maximal volume is based on the following result of the first author:

\begin{theorem}
\label{GuT}
(\cite {Gursky1}, {\em Theorem B}) If the Yamabe invariant $Y(M^4,[g]) \geq 0$, then 
\begin{align}
\label{Qineq}
\int_{M^4} \sigma_2(g^{-1}A) dvol = \int_{M^4} (-\frac{1}{8}|E|^2 + \frac{1}{96}R^2) dvol \leq 4\pi^2.
\end{align}
Furthermore, equality holds if, and only if, $(M^4,g)$ is conformally equivalent to the round
sphere.
\end{theorem}

To prove Theorem \ref{4dest},
suppose $g \in \Gamma_k^{+}(M^4)$ ($k\geq3$) satisfies
\begin{align*}
\sigma_k^{1/k}(g^{-1}A) \geq \sigma_k^{1/k}(S^4).
\end{align*}
From the Newton-Maclaurin inequality it follows that 
\begin{align*}
\sigma_2^{1/2}(g^{-1}A) \geq \sigma_2^{1/2}(S^4).
\end{align*}
Therefore,
\begin{align*}
\begin{split}
4\pi^2 &\geq \int_{M^4} \sigma_2(g^{-1}A) dvol \\
&\geq \sigma_2(S^4) vol(M^4,g) \\
&= \frac{3}{2} vol(M^4,g),
\end{split}
\end{align*}
and consequently $vol(M^4,g) \leq \frac{8}{3}\pi^2 = vol(S^4)$.  

Now suppose equality is attained in (\ref{sharp4}).  Then there is a sequence of metrics
$\{g_j\} \subset \Gamma_k^{+}(M^4)$ with $\sigma_k^{1/k}(g_j^{-1}A_{g_j}) \geq \sigma_k^{1/k}(S^4)$ and
$vol(M^4,g_j) \to vol(S^4) = \frac{8}{3}\pi^2$ as $j \to \infty$.  Therefore, $\sigma_2^{1/2}(g_j^{-1}A_{g_j}) \geq 
\sigma_2^{1/2}(S^4)$, and appealing once more to (\ref{Qineq}) we have
\begin{align*}
\begin{split}
4\pi^2 &\geq \int_{M^4} \sigma_2(g_j^{-1}A_{g_j}) dvol \\
&\geq \sigma_2(S^4) vol(M^4,g_j) \to 4\pi^2.
\end{split}
\end{align*}
Since $\int \sigma_2$ is conformally invariant, we conclude that
\begin{align*}
\int_{M^4} \sigma_2(g_j^{-1}A_{g_j}) dvol = 4\pi^2
\end{align*}
for each $j$.  By Theorem \ref{GuT}, $(M^4,g)$ is conformally equivalent to the
round sphere.

Theorem \ref{bpluscase} is a consequence of the following estimate:

\begin{theorem}
\label{GAn}
(\cite{Gursky2}, {\em Theorem 1})  If $M^4$ is a smooth, compact, orientable four-manifold with $b^{+} > 0$, then for any
metric $g$ of positive scalar curvature the Weyl tensor satisfies
\begin{align}
\label{WeylL2}
\int_{M^4} |W^{+}|^2 dvol \geq \frac{4}{3}\pi^2\left( 2\chi(M^4) + 3\tau(M^4) \right).
\end{align}
Furthermore, equality holds if, and only if $(M^4,g)$ is conformal to a positive K\"ahler-Einstein metric.   
In this case, $M^4$ is diffeomorphic to either $S^2 \times S^2$, $\mathbf{CP}^2$, or $\mathbf{CP}^2
\# m(-\mathbf{CP}^2)$ with $3\leq m \leq 8$.  
\end{theorem}

Suppose $b^{+} > 0$ and $g \in \Gamma_k^{+}(M^4)$ ($k\geq3$).  In particular, this implies that $g$ has positive scalar curvature
(i.e., $\sigma_1(g^{-1}A) > 0$).  
Combining the Chern-Gauss-Bonnet and signature formulas,
\begin{align}
\label{CGBS}
2\pi^2 \left( 2\chi(M^4) + 3\tau(M^4) \right) = \int_{M^4} |W^{+}|^2 dvol + 2 \int_{M^4} \sigma_2(g^{-1}A) dvol.
\end{align}   
Combining (\ref{CGBS}) and (\ref{WeylL2}), $g$ satisfies 
\begin{align*}
\frac{1}{3}\pi^2 \left( 2\chi(M^4) + 3\tau(M^4) \right) \geq \int_{M^4} \sigma_2(g^{-1}A) dvol.
\end{align*}
If we normalize $g$ so that 
\begin{align*}
\sigma_k^{1/k}(g^{-1}A) \geq \sigma_k^{1/k}(S^4),
\end{align*}
then the Newton-Maclaurin inequality implies  
\begin{align*}
\sigma_2^{1/2}(g^{-1}A) \geq \sigma_2^{1/2}(S^4).
\end{align*}
Therefore,
\begin{align*}
\begin{split}
\frac{1}{3}\pi^2 \left( 2\chi(M^4) + 3\tau(M^4) \right) & \geq \int_{M^4} \sigma_2(g^{-1}A) dvol \\
& \geq \sigma_2(S^4) vol(M^4,g) \\
&= \frac{3}{2} vol(M^4,g),
\end{split}
\end{align*}
and it follows that
\begin{align*}
\Lambda_k(M^4,[g]) \leq \frac{2}{9}\pi^2 \left( 2\chi(M^4) + 3\tau(M^4) \right).
\end{align*}
This proves Theorem \ref{bpluscase}

\bibliography{GV_AIMfinal_references}
\noindent
\small{\textsc{Department of Mathematics, University of Notre Dame, 
Notre Dame, IN 46556}}\\
{\em{E-mail Address:}} \ {\texttt{mgursky@nd.edu}}\\
\\
\small{\textsc{Department of Mathematics, Massachusetts Institute
of Technology, Cambridge, MA 02139}}\\
{\em{E-mail Address:}} \
{\texttt{jeffv@math.mit.edu}}
\end{document}